\newif\ifdebug
\debugfalse

\newif\ifpreliminary
\preliminarytrue

\documentclass[12pt,twoside]{article}

\usepackage{amsmath}
\usepackage{amstext}
\usepackage{amsbsy}
\usepackage{amssymb}
\usepackage{amsfonts}
\usepackage{latexsym}

\usepackage{graphicx}

\ifdebug
  \usepackage[notref,notcite]{showkeys}
\fi
\makeatletter
\@addtoreset{equation}{subsection}

\makeatother

\newcommand{\separate}{\medskip\noindent}

\def\newsection{ \separate
   \refstepcounter{subsection}
   {\large\bf \thesubsection\kern.3em}
}
\def\mytheorem#1{
   \separate{\large\bf Theorem#1:\kern.3em}}

\def\mylemma#1{
   \separate{\large\bf Lemma#1:\kern.3em}}

\def\mycorollary#1{
   \separate{\large\bf Corollary#1:\kern.3em}}

\def\myprop#1{
   \separate{\large\bf Proposition#1:\kern.3em}}

\def\myremark#1{
   \separate{\large\bf Remark#1:\kern.3em}}

\def\mydefinition#1{
   \separate{\large\bf Definition#1:\kern.3em}}

\def\myexample#1{
   \separate{\large\bf Example#1:\kern.3em}}

\def\Proof{\separate\underline{Proof:}\kern1em}

\newcommand{\field}[1]{\mathbb{#1}}
\newcommand{\C}{\field{C}}

\newcommand{\R}{\field{R}}

\newcommand{\N}{\field{N}}
\newcommand{\Z}{\field{Z}}

\newcommand{\X}{\mathbb X}

\newcommand{\MT}{\mbox{\small{$\widetilde{M}$}}}
\newcommand{\s}{\mbox{\boldmath $\sigma$}}

\newcommand{\CP}{\C \mathbb{P}^1}
\newcommand{\tr}{\mathrm{tr}}
\newcommand{\Id}{\mathrm{Id}}

\newcommand{\FID}{\mathcal{F}_{\mbox{\tiny{$\Id$}}}
(\MT)}

\def\QED{\hfill$\Box$}

\def\SL{{\mbox{\bf SL}(2,\C)}}

\def\GL{{\mbox{\bf GL}(2,\C)}}
\def\SU{{\mbox{\bf SU}(2)}}

\newcommand{\Utwo}{\mbox{\bf U}(2)}
\def\U1{{\mbox{\bf U}(1)}}

\def\sl{{\mbox{\bf sl}(2,\C)}}
\def\gl{{\mbox{\bf gl}(2,\C)}}
\def\su{{\mbox{\bf su}(2)}}

\newcommand{\n}{\noindent}

\def\LGC{{\Lambda_r \SL_\sigma}}

\def\LGP{{\Lambda_r^{\mbox{\tiny{$+$}}} \SL _\sigma}}

\def\LGU{{\Lambda_r \SU_\sigma}}

\def\Lhol{{\Lambda_{\mbox{\tiny{$\C^*$}}}\SU_\sigma}}

\def\LAC{{\Lambda_r \sl _\sigma}}

\def\Lsu{{\Lambda \su _\sigma}}
\def\LGe{{\Lambda_r \su _\sigma}}

\newcommand{\pot}{\Lambda \Omega(M)}
\newcommand{\pott}{\Lambda \Omega(\MT)}

\newcommand{\gauge}{\mathcal{G}_r (\MT)}

\def\s{{\mbox{\bf $\sigma$}}}

\def\l{{\mbox{\footnotesize \sc l}}}

\ifdebug
\else
\fi
\begin{document}

\begin{center}
{\LARGE Dressing preserving the fundamental group}

\vskip1cm
\begin{minipage}{6cm}
\begin{center}
Josef~Dorfmeister \\ 
TU M\"{u}nchen\\ 
Zentrum Mathematik \\ 
Boltzmannstr. 3 \\ 
D-85747 Garching, Germany.
\end{center}
\end{minipage}
\begin{minipage}{6cm}
\begin{center}
Martin~Kilia$\mbox{n}^*$ \\ 
Mathematical Sciences \\
University of Bath \\ 
Claverton Down \\
Bath, BA2 7AY, UK.
\end{center}
\end{minipage}
\vspace{0.5cm}

\end{center}
\footnotetext[0]{
  $^*$ Supported by CEC Contract 
HPRN-CT-2000-00101 EDGE and EPSRC GR/S28655/01. \date{\today}}

\begin{abstract}
In this note we consider the 
relationship between the dressing action and 
the holonomy representation in the context of 
constant mean curvature surfaces. 
We characterize dressing elements that preserve 
the topology of a surface and discuss dressing 
by simple factors as a means of adding bubbles 
to a class of non finite type cylinders.
\end{abstract}

\refstepcounter{subsection}

\textbf{Introduction.} 
The equation for a harmonic map from a Riemann 
surface to a Riemannian symmetric space has a 
zero-curvature representation, and so 
corresponds to a loop of flat connections. 
Uhlenbeck discovered in her study \cite{Uhl} of 
harmonic maps into a compact Lie group 
$G$ that such maps correspond to certain holomorphic 
maps into the based loop group of $G$ and used 
this to define the \textit{dressing action} of 
a certain loop group on the space of harmonic maps. 
Dressing, or the \emph{vesture method} \cite{ZakS2}, 
was first developed in soliton theory to generate 
new solutions from old by solving a matrix Riemann 
problem. Unfortunately, the new solution does not 
automatically inherit properties of the old, 
such as domain, periods or asymptotics nor is 
it easy to control such properties when solving 
a matrix Riemann problem. 

When the target is the two dimensional 
round sphere, harmonic maps are precisely 
the Gauss maps of surfaces with constant mean 
curvature \cite{RuhV}. In the article \cite{DorPW} 
a method was presented by which all conformally 
immersed surfaces with constant mean curvature 
(CMC) can be obtained. 
The construction involves solving a meromorphic 
linear differential system with values in 
a loop group and then Iwasawa decomposing the 
solution to obtain the extended unitary frame of 
the surface. Both these steps make it 
difficult to keep track of the topology. 
Variation of the initial condition 
corresponds to the dressing action and is 
an integral part of the theory. 
Not only can dressing be used to close periods, but 
also to generate new CMC surfaces from old ones 
without altering the topology. In this work, 
we look into the relationship 
between dressing and topology in the framework of 
meromorphic ODE's and loop group factorizations 
in the context of CMC surfaces.
 
The question, when a conformal CMC immersion 
factors through a given surface has been discussed 
for tori \cite{Bob:tor}, \cite{DorH:per}, 
\cite{PinS}, cylinders \cite{DorH:cyl}, 
\cite{Kil:thesis}, 
\cite{KilMS} and trinoids \cite{BobPS}, 
\cite{DorW:noids}, \cite{KilMS}, \cite{Sch:tri}. 
It will be interesting to see more results in 
this direction, with particular emphasis on 
complete, non-compact surfaces with genus and 
compact surfaces of genus $g \geq 2$.

Sterling \& Wente \cite{SteW} discovered that 
one can add ``bubbles'' to a round cylinder 
while preserving the mean curvature by 
Bianchi-- B\"{a}cklund transformations. 
In \cite{Kil:thesis} it was shown that these 
``bubbletons'' can be produced by dressing 
from the standard round cylinder by the so called 
{\em simple factors} of Terng and Uhlenbeck 
\cite{TerU}. Subsequently, it was shown in 
\cite{Kob} that the simple bubbletons in 
\cite{Kil:thesis} coincide with those of Sterling 
and Wente \cite{SteW}. Recently, A. Mahler 
\cite{Mah} has shown that indeed 
Bianchi -- B\"{a}cklund transformations  
can be achieved by dressing. Further, it was 
recently shown in \cite{KilSS} that it is 
also possible to dress CMC surfaces homeomorphic 
to the n-punctured sphere with suitably chosen 
simple factors without changing the topology.
 
Every CMC immersion of a Riemann surface $M$ induces
a 'monodromy representation' 
$\chi = \chi (\gamma,\,\lambda )$ of the 
fundamental group $\pi_1 (M)$ of $M$ with values 
in a unitary loop group. Moreover, as a function 
of the loop parameter $\lambda$, the monodromy 
matrices are holomorphic on $\C^*$. 
Dressing a given immersion with some arbitrary 
matrix will in general destroy the topology and 
the dressed surface will have a trivial 
fundamental group. Thus the question is: 
For which dressing matrices will the dressed 
immersion have the same fundamental group?

To explain this in more detail we note that in 
our method we associate with a given CMC immersion 
$f :M  \rightarrow \R^3$ its associated family 
$f_\lambda :\MT \to \R^3$, where $\MT$ 
denotes the universal cover of $M$ and 
$\lambda \in S^1$. This way we obtain the extended 
frame $F = F(z,\, \bar{z},\, \lambda )$ for 
$\lambda \in S^1$ and $z \in \MT$.
If $\gamma \in \pi_1(M)$ and we also denote the 
corresponding deck transformation by $\gamma$ and 
write 
$\gamma^* F = F(\gamma(z),\,\overline{\gamma(z)},\,\lambda)$  
then we have $\gamma^*F = \chi\,F\,k$ for some smooth 
$k:\MT \to \U1$. 
Since, by assumption, the immersion $f_1 = f$ 
descends to $M$, the monodromy matrices satisfy 
the two closing conditions 
$ \chi (\gamma,\, 1) = \pm \Id$ and 
$\left.\partial_{\lambda} \chi \right|_{\lambda = 1}
 = 0$ for all $\gamma \in \pi_1(M)$.

If $h = h(\lambda )$ is some dressing element 
and $\hat F$ denotes the dressed extended 
frame and we write $hF = \hat{F} V_+$ for 
the dressing equation, then $\gamma^* \hat{F} = 
h\,\chi\,h^{-1}\,\hat{F}\, V_+\,\gamma^* V_+^{-1}$. 
Thus $\hat{F}$ has some monodromy 
$\hat{\chi}$ if and 
only if $\hat{L} \hat{F} = \hat{F} \hat{W}_+$, 
where $\hat{W}_+ = 
V_+ \, \gamma^* V_+ ^{-1}$ and 
$\hat{L} = ( h \chi h^{-1})^{-1}
\hat{\chi}$. As a consequence, 
$\hat{\chi} = h \chi h^{-1} \hat{L} = 
h \chi L h^{-1}$, where $L = h^{-1} \hat{L} h$.

We prove 
that 
a dressing matrix $h$ preserves the fundamental 
group of an immersion if and only if for 
every $\gamma \in \pi_1 (M)$ there exists some 
matrix $L = L(\gamma,\,\lambda)$, such that 

(i) $LF = F W_+$,\vspace{-1mm}

(ii) $h \chi h^{-1}$ is holomorphic for 
  $\lambda \in \C^*$,\vspace{-1mm} 

(iii) $h \chi L h^{-1}$ is unitary on $S^1$, 
  and \vspace{-1mm}

(iv) The closing conditions are satisfied for 
  $h\chi L h^{-1}$.

Clearly, the situation simplifies considerably, 
if the original immersion has an umbilic point, 
since then $L = \pm \Id$ \cite{DorH:dre}. 
Let us assume this for now. 
Considering in this case condition (ii) we see 
that $\chi$ is holomorphic for 
$\lambda \in \mathbb{C}^*$ and also 
$h \chi h^{-1}$ needs to be holomorphic on 
$\mathbb{C}^*$, in spite of the fact
that $h$ may only be defined on some circle 
of radius $ 0 < r \leq 1$. 
We prove 
that in 
the situation discussed here we can assume that 
$h = \mathcal{M}\,\mathcal{C}$, where 
$[\,\mathcal{C}, \,\chi \,] = 0$ and $\mathcal{M}$ 
is meromorphic on some open dense subset 
$\mathbb{S} \subset \mathbb{C}^*$, 
which contains an open annulus about $S^1$.
Moreover, $\mathbb{S}$ is solely defined 
by the eigenvalues of $\chi$.
In view of condition (iii) we further show 
that actually 
$\mathcal{M}$ can be chosen so that it is 
unitary on $S^1$.

This gives altogether a fairly complete 
description for the case, where an umbilic 
point exists. Similarly complete is the case, 
dealt with in Theorem \ref{th:h-split_cyl}, 
where $\chi$ is the monodromy of the standard 
round cylinder. For the general umbilic free 
case we have, at this point, only some partial 
results.

Let us briefly outline the contents of this paper. 
The first chapter sets the scene by defining the 
relevant loop groups, dressing and gauge actions 
and recalls the DPW representation \cite{DorPW} 
of CMC surfaces. 
In the second chapter we consider how automorphisms 
affect maps at the various levels of the DPW 
construction as a prerequisite to understanding 
monodromy. In chapter 3 we present a factorization 
as a necessary condition on the dressing matrix 
to ensure that the dressed surface retains the 
topology of the original surface. 
In the fourth chapter we apply our 
methods to the dressing orbit of the vacuum 
and give a self contained account of twisted 
simple factors. We conclude this work by 
discussing how dressing with simple factors is 
related to the monodromy and apply these methods 
to a class of CMC cylinders with umbilics.

The second author wishes to thank F Burstall for 
many useful discussions and the Department of 
Mathematical Sciences at the University of Bath 
for its hospitality.

\section{Notation and basic results}\label{basics}
\message{[basics]} 
We begin by collecting some well known results 
on loop groups and the dressing action. 
We shall use the following notation 
for diagonal and off--diagonal matrices:
\begin{equation*}
\mathrm{diag}[u,\, v] = 
\bigl( \begin{smallmatrix} 
u & 0 \\ 0 & v \end{smallmatrix} \bigr),\,
\mathrm{off}[u,\, v] = 
\bigl( \begin{smallmatrix} 
0 & u \\ v & 0 \end{smallmatrix} \bigr). 
\end{equation*}
\newsection{\bf Loops.} \label{loopgroups}
For real $r \in (0,1]$, let 
$C_r = \{ \lambda \in \C : | \lambda | = r \}$ 
and denote the $r$--Loop group of $\SL$ by
\begin{equation*}
  \Lambda_r \SL = \left\{ g : C_r \to \SL
  \mbox{ smooth} \right\}.
\end{equation*}
We have an involution on maps $C_r \to \gl$ 
defined by 
\begin{equation} \label{eq:sigma}
  \sigma : g(\lambda) \mapsto 
  \sigma_3 g(-\lambda) \sigma_3^{-1},\,
  \sigma_3 = \mathrm{diag}[1,-1],
\end{equation}
and denote the \textit{twisted} $r$--Loop group 
of $\SL$ by
\begin{equation*}
  \LGC = \left\{ g \in \Lambda_r \SL: 
  \, \sigma g = g \right\}.
\end{equation*}
Analogously, one can define the Lie algebras of 
these groups, denoted by
$\LAC$. To make these loop groups 
complex Banach Lie groups, we equip them, 
as in \cite{DorPW},
with some $H^s$ topology for $s>1/2$ or some 
(possibly weighted) Wiener topology. 
Elements of these twisted loop
groups are matrices whose off--diagonal entries 
are odd functions, 
while the diagonal entries are even functions of 
the parameter $\lambda$. 
All entries are in the Banach algebra 
${\cal A}_r$ of
$H^s$--smooth functions or of finite 
(possibly weighted) Wiener norm on $C_r$.
Furthermore, we will use the following subgroups 
of $\LGC$: Let 
$I_r = \{ \lambda \in \C : |\lambda|<r \}$ 
and denote
\begin{equation*}
  \LGP = \left\{ g \in \LGC : g 
  \text{ extends analytically to } I_r \right\}.
\end{equation*}
Let $A_r = \{ \lambda \in \C : 
r<|\lambda|<1/r \}$, and by abuse of 
notation we denote 
\begin{equation*}
  \LGU = \left\{ g \in \LGC :  g 
  \text{ extends analytically to $A_r$ 
  and}\left. g \right|_{S^1} \in \SU \right\}.
\end{equation*}
Let $g:A_r \to \SL_\sigma$ be analytic. 
Then $g \in \LGU$ if and only if 
$\varrho \, g = g$ for 
\begin{equation} \label{eq:varrho}
  (\varrho g)(\lambda) :=
  \overline{g(1/\bar{\lambda})}^{t-1}.
\end{equation}
or alternatively, $g \in \LGU$ if and only if 
\begin{equation} \label{eq:unitarity_def}
  (g^*)(\lambda) := 
  \overline{g(1/\bar{\lambda})}^t = 
  g(\lambda)^{-1}.
\end{equation}    
For $r=1$ we will always omit the subscript '$r$'.
All the groups above are Banach Lie subgroups 
of $\LGC$. The corresponding Banach Lie 
sub algebras of $\LAC$ are defined analogously.


\newsection{\bf Iwasawa decomposition.} 
Of fundamental importance in our investigation 
is a certain Loop group factorization. 
We modify a result from \cite{McI} to obtain:

\separate  Multiplication 
$\LGU \times \LGP \to \LGC$
is a real analytic surjection. 
An associated splitting $g = FB$ 
of $g \in \LGC$ with $F \in \LGU$ and $B \in \LGP$ 
will be called an Iwasawa decomposition of $g$. 
Note that
\begin{equation} \label{eq:intersect}
  \LGU \cap \LGP = \U1.
\end{equation}
If we demand that $B(\lambda=0)$ has positive real 
diagonal entries, then the multiplication map 
above is a real analytic diffeomorphism onto.
We will call this decomposition the "unique" Iwasawa
decomposition.


\newsection{\bf Untwisted loops.} 
We shall be primarily working with the twisted 
loop groups and algebras. On occasion it is 
advantageous to switch to the untwisted setting 
via the isomorphism between untwisted and twisted 
$r$--loops: 
let $X_u \in \Lambda_r \SL$ be an 
untwisted loop. 
Then the corresponding twisted loop is given by 
$X_t(\lambda) = D\,X_u(\lambda^2)\,D^{-1}$ 
where $D=\mathrm{diag}[\,\sqrt{\lambda},\,
1/\sqrt{\lambda}\,]$. In matrix notation, this 
isomorphism is given by
\begin{equation} \label{eq:twist}
  \begin{pmatrix} a(\lambda) & b(\lambda) \\
  c(\lambda) & d(\lambda) \end{pmatrix} \cong 
  \begin{pmatrix} a(\lambda^2) & 
  \lambda\, b(\lambda^2) \\
  \lambda^{-1} c(\lambda^2) & 
  d(\lambda^2)\end{pmatrix}. 
\end{equation}
One technical issue in switching between untwisted 
and twisted loops is that under this twisting 
isomorphism, an untwisted positive loop is not 
automatically positive when twisted, since the 
untwisted loop at $\lambda = 0$ need not be 
upper triangular. To circumvent this issue, 
we must first perform a Gram-Schmidt 
orthonormalization at $\lambda = 0$ 
before twisting the loop. 

More precisely, consider an untwisted loop 
$\Phi_u \in \Lambda_r\SL$ with an $r$-Iwasawa 
decomposition  $\Phi_u = F_u\,B_u$ 
with untwisted $F_u \in \Lambda_r \SU$ and 
$B_u \in \Lambda_r^+ \SL$. Now consider the 
twisted loop $\Phi_t(\lambda) = D\,\Phi_u(\lambda^2)\,D^{-1}$ 
obtained from $\Phi_u$ via the 
isomorphism \eqref{eq:twist} and let 
$\Phi_t = F_t\,B_t$ be an $r$-Iwasawa decomposition 
with twisted $F_t \in \LGU$ and $B_t \in \LGP$. 
If $B_u$ has an expansion 
$B_u = B_0 + \lambda B_1 + \ldots$, 
then as $B_0 \in \SL$, we may write  
$B_0 = Q\,R$ via Gram-Schmidt orthonormalization 
with $Q \in \SU$ and $R \in \SL$ upper triangular. 
Consequently, $Q^{-1}B_u \in \Lambda_r^+\SL$ and 
$\left.Q^{-1}B_u \right|_{\lambda=0} = R$ and thus 
$D\,Q^{-1}B_u(\lambda^2)\,D^{-1} \in \LGP$. 
Furthermore, we then have that 
\begin{equation}
  \Phi_t(\lambda) = \left( D\,F_u(\lambda^2) 
  \,Q\,D^{-1} \right) \left( D\,Q^{-1} 
  B_u(\lambda^2) \,D^{-1} \right)
\end{equation}
with $D\,F_u(\lambda^2)\,Q\,D^{-1} \in \LGU$ and 
$D\,Q^{-1}B_u(\lambda^2)\,D^{-1} \in \LGP$ is an 
$r$-Iwasawa decomposition. It is in this sense that 
we may carry an Iwasawa decomposition in the 
untwisted setting over to the twisted setting.  


\newsection{\bf DPW method.} \label{sec:DPW}
Let $\Omega(M)$ denote the holomorphic 
$1$--forms on a Riemann surface $M$ and define 
the 
\begin{equation*}
  \pot = \Omega(M) \otimes \left\{ 
  \xi : \C^* \to  \sl 
  \mbox{ holomorphic}: \xi(\lambda) = 
  \sum_{j \geq -1} 
  \xi_j \lambda^j,\,\s \xi = \xi \right\}.
\end{equation*}
CMC surfaces come in $S^1$ families, 
the \textit {associated family}. 
The DPW representation \cite{DorPW} constructs 
all conformal CMC immersions of the universal 
cover $\MT$ in the following three steps: 
Let $\xi \in \pott$, $\tilde{z}_0 \in 
\MT$ and $\Phi_0 \in \LGC$ for some $r \in (0,1]$. 
To avoid totally umbilical surfaces, we assume 
$\det \xi_{-1} \not \equiv 0$.

\noindent 1. Solve the initial value problem 
\begin{equation}\label{eq:IVP}
  d\Phi = \Phi \xi,\, \Phi(\tilde{z}_0)= \Phi_0
\end{equation}
to obtain a unique \textit{holomorphic frame} 
$\Phi : \MT \to \LGC$.

\noindent 2. Iwasawa decompose the map 
$\Phi : \MT \to \LGC$ point-wise on $\MT$ 
\begin{equation} \label{eq:Iwasawa}
  \Phi(z,\,\lambda) = 
  F(z,\,\bar{z},\,\lambda) \, 
  B(z,\,\bar{z},\,\lambda) 
\end{equation}
We will always assume 
that the factors $F$ and $B$ in \eqref{eq:Iwasawa} 
are real analytic in $z$. 
The map $F : \MT \to \LGU$ will be called 
\textit{unitary frame}. If we use the unique 
Iwasawa decomposition, then the factors $F$ and $B$ 
are automatically real analytic.

\noindent 3. Let $H \in \R^*$ and 
$\partial_\lambda = 
\tfrac{\partial}{\partial \lambda}$. 
Plug $F$ into the Sym-Bobenko formula 
\begin{equation} \label{eq:Sym}
  f_\lambda = - \tfrac{1}{2H}\left( 
  i\lambda \tfrac{\partial F}{\partial\lambda} 
  \,F^{-1} + \tfrac{i}{2} F \sigma_3 F^{-1} \right)
\end{equation}
to obtain (possibly branched) conformal 
CMC immersions $f_\lambda :\MT \to \Lsu$ with mean 
curvature $H$, that is, for each 
$\lambda_0 \in S^1$ we have a conformal CMC 
immersion $f_{\lambda_0} :\MT \to \su \cong \R^3$. 
Note, $f_\lambda (z)$ is branched at $z_0$ if 
and only if for $\xi_{-1} = \mathrm{off}[a,\,b]$ we 
have $a(z_0)=0$.

\newsection{\bf Frames.} We call $\xi \in \pott$ a 
\textit{holomorphic potential}. 
If $\Phi_0 = \phi_u \phi_+$ is an Iwasawa 
decomposition, then it is not hard to see that 
the immersions obtained from 
$(\xi,\Phi_0,\tilde{z}_0)$ and 
$(\xi,\phi_+,\tilde{z}_0)$ 
differ by a $\lambda$-dependent rigid motion. 
Thus the choice of initial 
condition may be restricted to $\LGP$. 
Further, if $\Phi_0 \in \LGP$ 
then $F(\tilde{z}_0) \in \U1$ and 
$\Phi = F F(\tilde{z}_0)^{-1} F(\tilde{z}_0)B$. 
Hence we may assume without loss of generality that
\begin{equation} \label{eq:F_initial}
  F(\tilde{z}_0) \equiv \Id \mbox{ for all }
  \lambda \in S^1.
\end{equation}
Let $\sl = \mathfrak{k} \oplus \mathfrak{p}$ 
be the Cartan decomposition induced by the 
involution $\s$ defined in \eqref{eq:sigma}. 
The specific form of $\xi$ ensures that the 
Maurer--Cartan form of $F$ acquires the form
\begin{equation} \label{eq:MC}
  F^{-1}dF = 
  \lambda^{-1} \alpha_{\mathfrak{p}}^{\prime} + 
  \alpha_{\mathfrak{k}} + 
  \lambda \, \alpha_{\mathfrak{p}}^{\prime\prime}  
\end{equation}
where $\alpha_{\mathfrak{p}}^\prime$ and 
$\alpha_{\mathfrak{p}}^{\prime\prime}$ are 
$(1,\,0)$ respectively $(0,1)$-forms on $\MT$ 
with values in $\mathfrak{p}$ and 
$\alpha_{\mathfrak{k}}$ takes values in 
$\mathfrak{k}$. 
For each $p \in \MT$, 
$\lambda \mapsto F(p,\lambda)$ is 
holomorphic on $\C^*$. Maps 
\begin{equation}
  F:\MT \to \Lhol := 
  \bigcap_{r \in (0,1]}\LGU
\end{equation}
with the property \eqref{eq:MC} are called 
\textit{extended unitary frames} 
and denoted by $\mathcal{F}(\MT)$.
We call the extended unitary frames 
for which \eqref{eq:F_initial} holds the 
\textit{normalized extended unitary frames} 
and denote these by
\begin{equation}
  \FID = \left\{ 
  F \in \mathcal{F}(\MT) : 
  F(\tilde{z}_0) = \Id \mbox{ for some } 
  \tilde{z}_0 \in \MT \right\}.
\end{equation}
%


\newsection{\bf Gauge.} \label{Gauge}
The DPW representation $(\xi,\Phi_0,\tilde{z}_0) 
\mapsto \mathcal{F}(\MT)$ is a surjective map 
\cite{DorPW}. Injectivity fails since the 
\textit{Gauge group} 
\begin{equation}
  \gauge = \{ G : \MT
  \to \LGP \mbox{ holomorphic } \},
\end{equation}
acts by right multiplication on the fibers of 
this map. On the level of the potential, this 
\textit{gauge action} is given by 
\begin{equation}\label{eq:gauge}
  \xi.G = G^{-1} \xi \ G + G^{-1} dG.
\end{equation}
A computation shows that if $\Phi$ solves 
\eqref{eq:IVP} with triple 
$(\xi,\Phi_0 ,\tilde{z}_0 )$ and 
$G \in \gauge$ then the triples 
$(\xi,\Phi_0 ,\tilde{z}_0 )$ and 
$(\xi.G,\Phi_0 G(\tilde{z}_0),\tilde{z}_0 )$ 
induce the same CMC immersions, assuming $H \neq 0$.


\newsection{\bf Dressing.} \label{Dressing}
For $r \in (0,1]$, $h \in \LGC$ and 
$F_0 \in \mathcal{F}(\MT)$ let
\begin{equation}\label{eq:dress}
  hF_0 = F \, B 
\end{equation}
be a point-wise Iwasawa decomposition of $hF_0$ in 
$\LGC$, where every factor is analytic in $z$. Then 
$F \in \mathcal{F}(\MT)$.

Note that $F:\MT \to \LGU$ and 
$B:\MT \to \LGP$ are not uniquely determined, 
since for any smooth $U:\MT \to \U1$, we again 
have an Iwasawa decomposition given by 
$hF_0 = FU \,U^{-1}B$. We shall write 
\begin{equation}
  F \in [h\#F_0]
\end{equation}
to signify that $F:\MT \to \LGU$ satisfies 
\eqref{eq:dress}. We call $[h\#F_0]$ the 
\textit{dressing class} of $F_0$ under $h$ 
and say that $F \in [h\#F_0]$ was obtained 
by \textit{dressing} $F_0 \in \mathcal{F}(\MT)$ 
by $h \in \LGC$. 

To preserve the base point condition 
\eqref{eq:F_initial}, one needs to restrict 
to dressing with elements $h \in \LGP$.
 
\mylemma{}\label{th:initial} 
Let $h \in \LGC$ and $F_0 \in \FID$ be a normalized 
extended unitary frame. Then there exists a 
normalized extended unitary frame in $[h\#F_0]$ if 
and only if $h \in \LGP$.

\Proof A straightforward computation, 
see Proposition 2.9 of \cite{BurP:dre}, 
shows that the Maurer--Cartan 
form of any element in $[h\#F_0]$ 
is again of the form \eqref{eq:MC}. 
The issue is that for
$F_0 \in \FID$ we have 
$F_0(\tilde{z}_0) = \Id$, while for 
$F \in [h\#F_0]$ we have apriori only 
$F(\tilde{z}_0) \in \U1$. 

If $hF_0 = FB$ is an Iwasawa decomposition, 
then $F_0(\tilde{z}_0) = F(\tilde{z}_0) = \Id$, 
implies $h = B(\tilde{z}_0) \in \LGP$. 

Conversely, if $hF_0 = FB$ is an 
Iwasawa decomposition, then so is
\begin{equation}
  hF_0 = F F(\tilde{z}_0)^{-1} F(\tilde{z}_0)B.
\end{equation}
Hence $F F(\tilde{z}_0)^{-1}\in [h\#F_0]$ and 
satisfies \eqref{eq:F_initial}. Thus 
$F F(\tilde{z}_0)^{-1} \in \FID$. \QED

We thus have a left action of $\LGP$ 
on $\FID$ and shall write 
$F = h \# F_0$ to signify that $F \in [h\#F_0]$ 
with $F(\tilde{z}_0) = \Id$. Evaluating $hF_0 = 
h\#F_0 \, B$ at $\tilde{z}_0$ also gives 
$h=B(\tilde{z}_0)$. 
%

\newsection{\bf Isotropies.}
We denote the isotropy group of a map
$F: \MT \to \LGU$ under dressing by 
\begin{equation}
  \mathrm{Iso}_r(F) = \left\{ h \in \LGC : 
  F \in [h\# F] \right\}.
\end{equation} 
For a holomorphic frame $\Phi:\MT \to \LGC$ 
we define the isotropy under dressing by
\begin{equation}
  \mathrm{Iso}_r(\Phi) = \{ g \in \LGC : 
  g \, \Phi = \Phi \, G \mbox{ for some } G 
  \in \mathcal{G}_r(\MT) \}.
\end{equation}
If $g \in \mathrm{Iso}_r(\Phi)$ with  
$g\Phi = \Phi G$ and $\Phi = F\,B$ 
is an Iwasawa decomposition of $\Phi$, 
then $gF = FBGB^{-1}$ with $BGB^{-1} \in \LGP$. 
Hence $g \in \mathrm{Iso}_r(F)$. Conversely, if 
$h \in \mathrm{Iso}_r(F)$, and $hF=FH$ is an Iwasawa 
decomposition, then $h\Phi = \Phi B^{-1}HB$ 
and thus 
\begin{equation} \label{eq:isotropies}
  \mathrm{Iso}_r(\Phi) = \mathrm{Iso}_r(F).
\end{equation}
Further, using \eqref{eq:intersect} it is 
straightforward to verify that 
\begin{equation} \label{eq:iso_intersect}
  \mathrm{Iso}_r(F) \cap \LGU = \U1 \mbox{ for } 
  F \in \FID.
\end{equation}
Umbilic points occur naturally on compact 
CMC surfaces of genus $g \geq 2$ as well 
as on complete, open CMC surfaces with more than 
two ends. In this context, we quote 

\mytheorem{ \cite{DorH:dre}}
\label{th:isotropy}
If the associated family of 
$F \in \mathcal{F}(\MT)$ has umbilics, 
then $\mathrm{Iso}_r(F) = \{ \pm \Id \}$ for all 
$r \in (0,\,1]$. \QED

\section{Symmetries}
The notion of symmetry for CMC immersions has been 
discussed in the articles 
\cite{DorH:per,DorH:sym2}. 
It turns out that symmetry can be defined on 
various levels in the context of the DPW 
representation.
The problem of dealing with symmetries in the 
DPW representation stem mostly from the fact 
that symmetries on the potential level are, 
by and large, defined as coming from symmetries 
on the immersion level and have not yet been 
defined completely intrinsically on the potential 
level. Since symmetries associated with 
automorphisms $Aut(\MT)$ of the universal 
cover are central to this article, we recall 
the basic facts, retaining the notation of 
section \ref{sec:DPW}. 


\newsection{\bf Definitions.} \label{symmetries}
Let $f : \MT \to \Lsu$ be a CMC $H \neq 0$ immersion 
generated by a triple $( \xi, \Phi_0, \tilde{z}_0 )$ 
with $\Phi$ its holomorphic frame and $F$ its 
unitary frame. In the rest of 
this section we will always assume 
$\gamma \in Aut(\MT)$ and for
a map G with domain $\MT$ we shall write 
$\gamma^* G = G \circ \gamma$.  

\n {\bf{Immersion level}}: 
A triple 
$(\gamma,\,\X,\,T)$ with 
$\X \in \LGU$ and $T \in \LGe$ 
is called a symmetry of $f$ if and only if 
\begin{equation}\label{eq:sym_immersion}
  \gamma^* f = \X \, f \, \X^{-1} 
  + T.
\end{equation}

\n {\bf{Unitary frame level}}: 
A pair $(\gamma,\,\X)$ with 
$\X \in \LGU$ 
is called a symmetry of $F$ if and 
only if there exists a smooth map 
$k \in C^\infty(\MT,\U1)$ such that  
\begin{equation}\label{eq:sym_frame}
  \gamma^* F = \X \,F \,k.
\end{equation}

\n {\bf{Holomorphic frame level}}: 
A pair $(\gamma,\,\X)$ with 
$\X \in \LGU$ 
is a symmetry of $\Phi$ if and only if 
there is a map $H \in \gauge$ such that
\begin{equation}\label{eq:sym_curve}
  \gamma^* \Phi = \X \,\Phi \,H.
\end{equation}

\n {\bf{Potential level}}: 
A pair $(\gamma,\,G)$ with 
$G \in \gauge$ is called a 
symmetry of $\xi$ if and only if 
\begin{equation}\label{eq:sym_potential}
  \gamma^* \xi = \xi.G \,.
\end{equation}
where $\xi.G$ denotes the gauge 
transformation of $\xi$ by $G$ as defined in 
\eqref{eq:gauge}.
á

\mytheorem{}\label{th:symmetries} 
Let $f : \MT \to \Lsu$ be a CMC $H\neq 0$ immersion 
generated by the triple 
$(\xi,\Phi_0,\tilde{z}_0 )$ with $\Phi$ and 
$F$ the holomorphic respectively unitary frame. 
\begin{enumerate}
  \item If $(\gamma,\,\X_1,\,T_1)$ and 
  $(\gamma,\,\X_2,\,T_2)$ are symmetries 
  of $f$, then $\X_2 = \pm \X_1$ and 
  $T_1 = T_2$. 
  \item If $(\gamma,\,\X_1)$ and $(\gamma,\,\X_2)$ 
  are symmetries of $F$ then $\X_2 = \pm \X_1$. 
  \item If $(\gamma,\,\X_1)$ and $(\gamma,\,\X_2)$ 
  are symmetries of $\Phi$ then $\X_2 = \pm \X_1$. 
\end{enumerate}
\Proof (i) If $(\gamma,\,\X_1,\,T_1)$ and 
$(\gamma,\,\X_2,\,T_2)$ are symmetries, then 
$\gamma^* f = \X_1 f \X_1^{-1} + T_1 = 
\X_2 f \X_2^{-1} + T_2$ implies 
$\X \, f \, \X^{-1} + T = f$ with 
$\X = \X_2^{-1}\X_1$ and 
$T = \X_2^{-1}(T_1-T_2)\X_2$. Since a CMC 
$H \neq 0$ surface in $\R^3$ is never contained 
in an affine plane, 
$\X \, f \, \X^{-1} + T = f$ implies 
$\X = \pm \Id$ and $T=0$. 

(ii) If $(\gamma,\,\X_1)$ and $(\gamma,\,\X_2)$ 
are symmetries of $F$, then 
$\gamma^* F = \X_1 F k_1 = \X_2 F k_2$. 
Plugging $\X_1 F k_1$ and $\X_2 F k_2$ into 
the Sym--Bobenko formula and equating yields
\begin{equation}
  \X_1 f \X_1^{-1} - \tfrac{i\lambda}{2H}
  (\partial_{\lambda} \X_1 ) \X_1^{-1} = 
  \X_2 f \X_2^{-1} - \tfrac{i\lambda}{2H}
  (\partial_{\lambda} \X_2 ) \X_2^{-1}.
\end{equation}
Now $\X_2^{-1} \X_1 \,f \, \X_1^{-1} \X_2 + 
\tfrac{i\lambda}{2H} \X_2^{-1} 
\left( (\partial_{\lambda} \X_2) \X_2^{-1} - 
(\partial_{\lambda} \X_1 ) \X_1^{-1} \right) \X_2 
= f$ and part (i) implies $\X_2 = \pm \X_1$.

(iii) If $(\gamma,\,\X_1)$ and 
$(\gamma,\,\X_2)$ are symmetries of $\Phi$ then 
$\gamma^* \Phi = \X_1 \Phi H = \X_2 \Phi G$ 
implies $\X_2^{-1} \X_1 \Phi = \Phi G H^{-1}$ 
and shows that 
$\X_2^{-1} \X_1 \in \mathrm{Iso}(\Phi)$. 
Let $\Phi = F B$ be an Iwasawa decomposition of 
$\Phi$. Then 
$\X_2^{-1} \X_1 \in \mathrm{Iso}(F)$, since 
$\mathrm{Iso}(\Phi) = \mathrm{Iso}(F)$ by 
\eqref{eq:isotropies}. 
Hence $\X_2^{-1} \X_1 = \pm \Id$ 
by \eqref{eq:iso_intersect}. \QED

For symmetries $(\gamma,\,\X_\gamma)$ and 
$(\mu,\,\X_\mu)$ of $F \in \mathcal{F}(\MT)$ and 
corresponding maps 
$k_\gamma,\,k_\mu$,  
we have $\mu^*\,\gamma^*\,F = 
\X_\gamma \, \X_\mu \,F\,k_\mu
  \,\mu^* k_\gamma$ and if we assume that 
$((\gamma \circ \mu),\,\X_{\gamma\mu})$ 
is a symmetry, then part (ii) of Theorem 
\ref{th:symmetries} implies $\X_{\gamma \mu} = 
\pm \X_\gamma \, \X_\mu$ and consequently
\begin{equation}
    k_{\gamma \mu} = \pm k_\mu \, \mu^* k_\gamma.
\end{equation}
%


\newsection{\bf Symmetries of $\Phi$.} 
We briefly investigate how symmetries 
of a holomorphic potential descend to symmetries 
of the corresponding holomorphic frame and note a 
simple consequence in case the potential is 
invariant. In view of these results and our 
subsequent inquiry into the relationship with the 
dressing action, we allow symmetries on the level 
of holomorphic frames of the form $(\gamma, X)$ 
with $X \in \LGC$.

\myprop{} Let $\gamma \in Aut(\MT)$ and 
$(\gamma,\,G)$ be a symmetry of 
$\xi \in \pott$. Let $\Phi:\MT \to \LGC$ 
be a solution to $d\Phi = \Phi \xi$. 
Then there exists $X \in \LGC$ 
such that $\gamma^* \Phi = X \Phi G$.

\Proof Since $\gamma^* \Phi$ and 
$\Phi G$ both solve the differential equation 
$dY = Y\gamma^*\xi$, there exists 
$X \in \LGC$ with $\gamma^*\Phi = X \Phi G$. 
More explicitly, if $\Phi$ solves 
$d\Phi = \Phi \xi,
\, \Phi(\tilde{z}_0) = \Phi_0$, then 
$X = \Phi(\gamma(\tilde{z}_0)) 
G(\tilde{z}_0)^{-1} \Phi_0^{-1}$.
\QED 

We shall again call a pair $(\gamma,\,X)$ with 
$\gamma \in Aut(\MT)$ and $X \in \LGC$ a 
symmetry of a map $\Phi:\MT \to \LGC$ 
if there exists a $G \in \gauge$ such that 
$\gamma^* \Phi = X\,\Phi \,G$.

\mycorollary{} \label{th:monophi} 
Let $\gamma \in Aut(\MT)$ and 
$\xi\in \pott$ with $\gamma^* \xi = \xi$. 
Let $\Phi$ be a solution to $d\Phi = \Phi \xi$. 
Then there exists $X \in \LGC$ such that 
$\gamma^* \Phi = X \Phi$. \QED


\newsection{\bf Symmetries of $F$.} 
In view of the previous section, 
we now characterize how symmetries 
of the holomorphic frame descend to symmetries 
of the corresponding extended unitary frame 
and show that for open Riemann surfaces, the 
co-cycle factor in equation \eqref{eq:sym_frame} 
can be gauged away.

\mylemma{}\label{th:symF} 
Let $(\gamma,\,X)$ be a symmetry 
of $\Phi:\MT \to \LGC$. Let 
$\Phi = F\,B$ be an Iwasawa decomposition 
of $\Phi$. Then there exists an $\X \in \LGU$ 
such that $(\gamma,\X)$ is a symmetry of 
$F$ if and only if there exists an element 
$L \in \mathrm{Iso}_r(\Phi)$ with $LX^{-1} \in \LGU$.

\Proof Let $G \in \gauge$ such that 
$\gamma^* \Phi = X \Phi G$. If $(\gamma,\X)$ 
with $\X \in \LGU$ is a symmetry of $F$, then by 
\eqref{eq:sym_frame} there 
exists a differentiable map $k:\MT \to \U1$ 
such that $\gamma^* F = \X F k$. In combination 
with $F = \Phi B^{-1}$, after rearranging, 
we obtain 
\begin{equation}
  \X^{-1} X \, \Phi = 
  \Phi B^{-1}k\gamma^*B \, G^{-1}.
\end{equation}
Then $L:= \X^{-1} X \in \mathrm{Iso}_r(\Phi)$ as 
$B^{-1}k\gamma^*B \, G^{-1} \in \gauge$ and 
$L \, X^{-1} = \X^{-1} \in \LGU$.

Conversely, let $L \in \mathrm{Iso}_r(\Phi)$ and 
$L \, X^{-1} \in \LGU$. Then there exists 
an element $H \in \gauge$ with 
$L\, \Phi = \Phi H$. By \eqref{eq:isotropies}, 
$L \in \mathrm{Iso}_r(F)$ and $LF=FV$ for 
$V = BHB^{-1} \in \gauge$. Then 
$\gamma^* \Phi = X \Phi G$ yields 
\begin{equation}
  \gamma^* F = X \,F\, B\, G\,  \gamma^* B^{-1} = 
  X\, L^{-1}F\, V\, B \,G\, \gamma^*B^{-1}.
\end{equation}
Define $\X:=X\,L^{-1}$ and 
$k:=VBG\gamma^*B^{-1}$. 
On the one hand, $k \in \gauge$ 
while on the other hand 
$k= F^{-1}\X^{-1}\gamma^*F$ 
takes values in $\LGU$, since 
$\X \in \LGU$ by assumption. From equation 
\eqref{eq:intersect} we conclude that 
$k:\MT \to \U1$. \QED

\mytheorem{}
Let $M$ be an open Riemann surface 
with Fuchsian group $\Gamma$ and let 
$f_\lambda: \MT \to \R^3$ be an associated family 
of CMC $H \neq 0$ immersions. 
Then there exists an extended unitary frame 
$F \in \mathcal{F}(\MT)$ for $f_\lambda$ such 
that for every $\gamma \in \Gamma$ there exists 
$\X \in \Lhol$ such that $\gamma^*F = \X \,F$.

\Proof First we apply \cite{DorH:sym2}, 
Theorem 2.3 and infer that there 
exists some extended frame 
$\tilde{F} \in \mathcal{F}(\MT)$ for $f_\lambda$ 
such that for every $\gamma \in \Gamma$ we have 
$\gamma^* \tilde F =
\X \,\tilde{F} \,k$, with $\X \in \Lhol$ and 
$k:\MT \to \U1$ as in \eqref{eq:sym_frame}. 
Writing $\tilde{F} = \Phi \,B$ with 
$\Phi:\MT \to \LGC$ holomorphic and 
$B \in \gauge$, we obtain 
$\gamma^* \Phi = \X\,\Phi\,W$ 
with $W \in \gauge$. It is easy to see that 
$W$ satisfies the co-cycle condition. Therefore, 
since the open Riemann surface $\MT$ is Stein, 
by \cite{Bun}, $W$ is a co-boundary and whence 
$W = W\,\gamma^*W^{-1}$. Hence 
$\tilde{\Phi} = \Phi \,W$ satisfies 
$\gamma^* \tilde{\Phi} = \X \, \tilde{\Phi}$. 
Iwasawa splitting $\tilde{\Phi} = F \hat{B}$ with 
$\hat{B}$ normalized such that the $\lambda^0$ 
coefficient in $\hat{B}$ has positive real entries 
we obtain $\gamma^* F = \X\,F$. \QED


\newsection{{\bf Symmetries \& Dressing.} 
\label{sec:sym} 
Next, we characterize how symmetries between two 
dressing equivalent unitary frames are related.

\mytheorem{} \label{th:sym1}
Let $F_0 \in \mathcal{F}(\MT)$ have symmetry 
$(\gamma,\X_0)$ and $h \in \LGC$. Then 
$F \in [h \# F_0]$ has symmetry $(\gamma,\X)$ 
if and only if there exists a 
$L \in \mathrm{Iso}_r(F)$ such that 
 $\X = h \X_0 h^{-1} L^{-1}$. 
Furthermore, if both $F_0,\,F \in \FID$, then 
$L \in \LGP$.

\Proof
Writing $hF_0 = F B$ for an Iwasawa decomposition 
and using $\gamma^* F_0 = \X_0 F_0 k_0$, we obtain 
$h \X_0 h^{-1} F  = 
  \gamma^* F \gamma^* B k_0^{-1} B^{-1}$. 
If $(\gamma,\X)$ is a symmetry of $F$ then 
$\gamma^* F = \X F k$, whence
\begin{equation} \label{eq:iso_con}
  \X^{-1} h \X_0 h^{-1} F =  
  F k \gamma^*B k_0^{-1} B^{-1}.
\end{equation}
Set $L :=\X^{-1} h \X_0 h^{-1}$ and 
$H := k \gamma^*B k_0^{-1} B^{-1}$. Clearly 
$L \in \LGC$ and $H \in \gauge$, so equation 
\eqref{eq:iso_con} reads $LF=FH$ and is an 
Iwasawa decomposition of $LF$. Hence 
$L \in \mathrm{Iso}_r(F)$, and by construction    
$\X = h \X_0 h^{-1} L^{-1}$.

Conversely, assume there exists a map 
$L \in \mathrm{Iso}_r(F)$ 
such that $L F = F H$ for some $H \in \gauge$ 
and $h \X_0 h^{-1} L^{-1} \in \Lhol$. 
Defining $\X = h \X_0 h^{-1} L^{-1}$ and 
$k = H B k_0 \gamma^* B ^{-1}$, a computation 
yields $\gamma^* F = \X F k$. 
A priori, $k:\MT \to \LGP$, but as 
$\gamma^* F,\, \X$ and $F$ take values 
in $\Lhol$ so does $k$. 
Hence $k:\MT \to \U1$ and consequently 
$(\gamma,\X)$ is a symmetry of $F$.
\QED

\mycorollary{}
If in addition to the assumptions of Theorem 
\ref{th:sym1}, $\mathrm{Iso}_r(F) = \{ \pm \Id \}$, 
then $(\gamma,\X)$ is a symmetry of 
$F$ if and only if $\X = \pm h\, \X_0 h^{-1}$.

\Proof
The assumption $\mathrm{Iso}_r(F) = 
\{ \pm \Id \}$ implies that for the maps $L,\,H$ in the 
proof of Theorem \ref{th:sym1} we have  
$L = \Id$ and $H \equiv \Id$, yielding the claim.
\QED


\newsection{\bf Groups of symmetries.}
Extending the results of section \ref{sec:sym} 
to groups $\Gamma \subset Aut(\MT)$ of symmetries 
one obtains:

\mytheorem{} \label{th:mainresult}
Let $F_0 \in \mathcal{F}_r(\MT)$ admit symmetries 
$(\gamma,\X_0(\gamma))$ for all 
$\gamma \in \Gamma$ and let $h \in \LGC$. 
Then $F \in [h \# F_0]$ 
has symmetries $(\gamma,\X(\gamma))$ for all 
$\gamma \in \Gamma$ if and only if there exist 
maps $L(\gamma) \in \mathrm{Iso}_r(F)$ such that 
$\X(\gamma)=h\X_0(\gamma) h^{-1} L(\gamma)^{-1}$. 
Moreover, if $\mathrm{Iso}_r(F) = \{ \pm \Id \}$, 
then $F$ admits symmetries $(\gamma,\X(\gamma))$ 
for all $\gamma \in \Gamma$ if and only if 
$\X(\gamma) = h \X_0(\gamma) h^{-1}$ for all 
$\gamma \in \Gamma$. \QED


\newsection{\bf Symmetries \& Isotropy} 
Let $(\gamma,\X)$ be a symmetry of 
$F \in \mathcal{F}(\MT)$ and 
$\gamma^* F = \X F k$. 
If $L \in \mathrm{Iso}_r(F)$ and $LF = FB$ for an 
Iwasawa decomposition, then 
$\X^{-1}L \X F = F k\gamma^*Bk^{-1}$. 
As $(k\gamma^*Bk^{-1})$ is $\LGP$--valued, 
$\X^{-1} L \X \in \mathrm{Iso}_r(F)$ 
and every symmetry 
$(\gamma,\X)$ induces an inner automorphism 
$S:\mathrm{Iso}_r(F) \to \mathrm{Iso}_r(F),\,
L \mapsto \X^{-1} L \X$. 

Now consider $F_0 \in \FID,\,
h \in \LGP$ and let $F = h \# F_0$. Assume that 
$(\gamma,\,\X_0)$ and $(\gamma,\,\X)$ are 
symmetries of $F_0$ respectively $F$. 
If further $h \in \mathrm{Iso}_r(F)$, then we 
obtain a map $S : \mathrm{Iso}_r(F) \to 
\mathrm{Iso}_r(F),\,h \mapsto 
\X^{-1} h \X_0 h^{-1}$ with the property 
$S(h_1 h_2) = S(h_1) h_1 S(h_2) h_1^{-1}$.


\section{Invariant potentials \& Monodromy}
We start our investigation of the relationship 
between monodromy and dressing by defining the 
notion of monodromy on the level of 
holomorphic and unitary frames in the DPW framework. 

\mydefinition{} \label{def:monodromies} 
Let $M$ be a connected Riemann surface 
with universal cover $\MT$ and $\Delta$ 
the group of deck transformations. 
Let $\xi \in \pott$ and $\Phi: \MT \to \LGC$ be a 
solution of $d \Phi = \Phi \xi$ and 
$\tau \in \Delta$. A loop 
$\varrho(\tau) \in \LGC$ such that 
$(\tau,\,\varrho(\tau))$ is a symmetry of 
$\Phi$ will be called a \textit{monodromy} of 
$\Phi$ with respect to $\tau$. 
Let $\Phi = F\,B$ be an Iwasawa decomposition. 
A loop $\chi(\tau) \in \LGU$ such that 
$(\tau,\,\chi(\tau))$ is a symmetry of 
$F$ is called a \textit{monodromy} of $F$ 
with respect to $\tau$.

\myremark{} A monodromy of a normalized extended 
unitary frame $F \in \FID$ takes values in $\Lhol$.
%
%

\newsection{\bf Existence of monodromies.} 
Apriori, it is not clear that 
for some $\tau \in \Delta$ there exist monodromies 
$\varrho(\tau)$ and 
$\chi(\tau)$ as in the definitions above. 
In case the 
corresponding potentials are invariant under 
$\Delta$ we can evoke Corollary \ref{th:monophi} 
and obtain the following sufficient conditions 
for the existence of monodromies.

\myprop{} (i) Let $\xi \in \pott$ and 
$\Phi$ be a solution of $d\Phi = \Phi\,\xi$ with 
initial condition $\Phi_0 \in \LGC$. 
Let $\tau \in \Delta$ with $\tau^*\xi = \xi$. 
Then there exists a monodromy 
$\varrho(\tau) \in \LGC$ such that 
$\tau^* \Phi = \varrho(\tau) \,\Phi$.

(ii) Let $F:\MT \to \LGU$ and $\alpha = F^{-1}dF$. 
If $\tau^*\alpha = \alpha$ for $\tau \in \Delta$, 
then there exists a monodromy $\chi(\tau) \in \LGU$ 
such that $\tau^* F = \chi(\tau) \,F$. \QED

It is shown in \cite{DorH:cyl} that CMC 
immersions of open Riemann surfaces $M$ 
can always be generated by $\Delta$--invariant 
potentials $\xi \in \pott$. 
If $\Phi = F B$ is the Iwasawa 
decomposition of a holomorphic frame 
$\Phi : \MT \to \LGC$, then we shall need 
to study the monodromy of $F$ in order to control 
the periodicity of the 
resulting immersion \eqref{eq:Sym}. 
Even if $F$ is obtained from a potential that 
is invariant under $\Delta$, we are a priori 
not assured that there exist loops 
$\chi(\tau) \in \LGU$ and maps 
$k(\tau):\MT \to \U1$ such that 
equation \eqref{eq:sym_frame} holds for all 
$\tau \in \Delta$ or even more strongly, 
that for $\alpha = F^{-1}dF$ we have 
$\tau^* \alpha = \alpha$ for all 
$\tau \in \Delta$. In analogy to Theorem 
\ref{th:sym1}, this issue is characterized 
by the following 

\mylemma{} \label{th:iso_problem} 
Let $(\xi, \Phi_0, \tilde{z}_0 )$ generate the 
holomorphic frame $\Phi:\MT \to \LGC$ and 
extended unitary frame $F \in \mathcal{F}(\MT)$. 
Let $\tau \in \Delta$ and $\varrho(\tau)$ 
be a monodromy of $\Phi$. 
Then the following are equivalent.
\begin{enumerate} 
  \item There exists a monodromy 
  $\chi(\tau) \in \LGU$ of $F$ with 
  respect to $\tau$.
  \item There exists an element  
  $L(\tau) \in \mathrm{Iso}_r(\Phi)$ 
  with $L(\tau)\varrho^{-1}(\tau) \in \LGU$. 
\end{enumerate}

\Proof We first show that (i) implies (ii). 
  Let $\chi(\tau) \in \LGU$ 
  and $k(\tau):\MT \to \U1$ such that 
  $\tau^* F = \chi(\tau) F\, k(\tau)$. 
  Using $F = \Phi B^{-1}$ and 
  $\tau^* \Phi = \varrho(\tau) \Phi$ and 
  rearranging, we obtain
  $\chi(\tau)^{-1}\varrho(\tau) \, \Phi = 
  \Phi\, B^{-1} k(\tau) \,\tau^* B$.
  Then $L(\tau) :=  \chi(\tau)^{-1}\varrho(\tau) 
  \in \mathrm{Iso}_r(\Phi)$ and 
  $L(\tau)\varrho(\tau)^{-1} \in \LGU$. 

Conversely, let $L(\tau) \in \mathrm{Iso}_r(\Phi)$ 
and assume $L(\tau)\varrho(\tau)^{-1} \in \LGU$. 
Then there exists an element 
$G(\tau) \in \gauge$ with 
$L(\tau) \, \Phi = \Phi \, G(\tau)$.
Multiplying this last equation on the right by 
$\tau^* \Phi^{-1}$ and rearranging gives
\begin{equation} \label{eq:iso_phi2} 
  F^{-1}L(\tau) \, \varrho(\tau)^{-1} \,\tau^*F 
  = B \, G(\tau) \, \tau^*B^{-1}.
\end{equation}
The left respectively right hand side of 
\eqref{eq:iso_phi2} takes values in $\LGU$ 
respectively $\LGP$. 
Hence, by \eqref{eq:intersect}, 
both sides are $\lambda$--independent and 
$\U1$--valued. Set $k(\tau):= F^{-1}L(\tau)\,
\varrho(\tau)^{-1}\tau^*F$ 
and $\chi(\tau) := \varrho(\tau)\,L(\tau)^{-1}$. 
Then $\tau^* F = \chi(\tau)\, F\, k(\tau)$ and 
$\chi(\tau) \in \LGU$. 
\QED

\mycorollary{} \label{th:chi=rho} 
Let $(\xi, \Phi_0, \tilde{z}_0 )$ generate the 
holomorphic frame $\Phi:\MT \to \LGC$ and 
extended unitary frame $F \in \mathcal{F}(\MT)$. 
Let $\tau \in \Delta$ and $\varrho(\tau)$ 
be a monodromy of $\Phi$ 
and $\chi(\tau)$ a monodromy of $F$ with 
respect to $\tau \in \Delta$. 
If $\mathrm{Iso}_r(F) = \{ \pm \Id \}$, 
then $\chi(\tau) = \pm \varrho(\tau)$. \QED
 

\newsection{\bf A factorization Theorem.} 
In this section we work exclusively with 
untwisted loops.Consider an analytic map $H:\C^* \to \SL$, 
for which $\tr H$ is not independent of $\lambda$.
Its eigenvalues are 
\begin{equation}
  \mu_{\pm} = \frac{1}{2} \left( \tr H \pm 
  \sqrt{(\tr H)^2 - 4}\right).
\end{equation}
We will need to have $\mu_{\pm}$ holomorphic on 
some open and dense subset $\mathbb S$ of $\C^*$.
To this end we slightly generalize the procedure 
used in \cite{DorH:cyl}. 

\myprop{}
Let $g:\C^* \to \C,\,g \not\equiv 0$ be holomorphic. 
Then there exists some connected, open, 
and dense subset $\mathbb S \subset \C^*$ 
such that there exists a well defined square 
root function $\sqrt{g}$ of $g$ on $\mathbb S$.

\Proof
Since $g(\lambda)$ does not vanish 
identically it has at most finitely 
many roots on $S^1$. Moving slightly inside, 
if necessary, we can choose any $0 < r \leq 1$ 
such that $g$ does not vanish on $C_r$.
We form a (scalar) Birkhoff splitting 
of $g$ on $C_r$: $g = g_- \lambda^N g_+$.
It is easy to verify that, since $g$ is 
holomorphic on $\C^*$, also $g_-$ and 
$g_+$ are holomorphic on $\C^*$. 
Now we introduce some cuts: Consider 
the roots of $g_+$. They are all contained 
in the complement of $I_r$.
We cut the Riemann sphere from every odd 
ordered root of $g_+$ to infinity. 
Similarly we cut the Riemann sphere from 
each odd ordered root of $g_-$ to the origin.
Thus $g_-$ and $g_+$ have well defined 
root functions on the cut domain.
If $N$ is odd, then we also need to 
introduce one further cut from the origin to the 
point at infinity. This resulting domain will 
be denoted by $\mathbb S$. 
\QED

Applying this to the function 
$g(\lambda) = {\tr H}^2 -4$ we obtain

\mycorollary{} \label{th:mu}
The eigenvalue functions $\mu_{\pm}(\lambda)$ 
of a loop $H \in \Lambda_{\C^*}\SL$ are analytic 
on a connected, open and dense set 
$\mathbb S \subset \C^*$, obtained by making 
radial branch cuts from odd ordered roots $\alpha$ 
of ${\tr H}^2 -4$ to the origin 
respectively the point at infinity, depending 
on whether $|\alpha| <r$ respectively 
$|\alpha| >r $. In addition there is 
possibly a cut from $0$ to $\infty$.
\QED

For our approach it will be convenient 
to diagonalize $H$
\begin{equation} \label{eq:diag}
  Y\, H \, Y^{-1}= \mathrm{diag}[\mu_+,\,\mu_-] 
\end{equation}
If the matrix entries $H_{1,2}$ and $H_{2,1}$ 
of $H$ do not vanish identically, 
then one can choose as diagonalizing matrix $Y$ 
( see  e.g. (3.5.18) in \cite{DorH:cyl})
\begin{equation}
  Y = \begin{pmatrix}
  1 &  -\frac{H_{12}}{u-iv-H_{11}}\\
  -\frac{H_{21}}{2iv}& \frac{u + iv-H_{22}}{2iv}
  \end{pmatrix},
\end{equation}
where  $H_{ij}$ denotes the $(i,j)-$coefficient 
of $H$ and $ 2u = \tr H$ and $ v^2 = u^2 -1$.
If one of the off-diagonal coefficients vanishes 
identically, then we choose the corresponding 
matrix $Y$ listed in \cite{DorH:cyl}. In this 
case the diagonal entries are $\mu_{\pm}$. Note 
that this works in our setting.(Not only in the 
somewhat more special setting of \cite{DorH:cyl}.)

From the specific form of the matrices 
$Y$ mentioned above we infer that $Y$ 
is meromorphic on $\mathbb S $, provided
the trace $\tr H$ of $H$ is not independent 
of $\lambda$.

The next result is crucial for our 
characterization of dressing matrices
preserving the fundamental group. As before, 
we denote an open annulus about $S^1$ by 
$A_r = \{ \lambda \in \C : r<|\lambda|<1/r \}$ and 
by $C_r = \{ \lambda \in \C : |\lambda| = r \}$ 
a circle of radius $r$. With the above notations 
we have the following

\mytheorem{} \label{th:h-split}
Let $h \in \Lambda_r^+\SL$ for some $0<r\leq1$ and 
$H  \in \Lambda_{\C^*}\SL$ such that the trace 
of $H$ is not constant. 
Let $\mathbb S$ be as in Corollary \ref{th:mu}. 
If $h\, H  \,h^{-1}$ is meromorphic on 
$\mathbb S \cap A_s$ for some 
$0 < s < r$, then $h$ can be factored, on a possibly 
segmented circle $C_l \cap \mathbb S$ for 
$s\leq l \leq r$, into 
\begin{equation}\label{eq:h-split}
  h = \mathcal{M}\, \mathcal{C},
\end{equation}
where $\mathcal{C}$ is  $\SL$--valued 
on $\mathbb S\cap C_l$ and 
$[\, \mathcal{C},\, H  \,] = 0$ there. Moreover, 
$\mathcal{M}$ is $\SL$--valued and meromorphic 
on the connected open set $\mathbb S \cap A_s$. 

\Proof From the above discussion, 
there exists some $Y$  
such that \eqref{eq:diag} 
holds. Let us write 
\begin{equation}
  h Y^{-1} = \begin{pmatrix} 
  a & b \\ 
  c & d \end{pmatrix}.
\end{equation}
Computing $h\,H \,h^{-1}$ and evaluating 
the diagonal entries gives that both 
\begin{align}
  &ad\mu_- -bc\mu_+ , \label{eq:r22}\\
  &ad\mu_+ - bc\mu_- \label{eq:r11}
\end{align}
are meromorphic on $\mathbb S \cap A_s$. 
Multiplying \eqref{eq:r22} by $\mu_+$ and 
\eqref{eq:r11} by $\mu_-$ and subtracting 
the resulting equations, 
implies that $bc(\mu_+^2 - \mu_-^2)$ 
is meromorphic on $\mathbb S \cap A_s$. 
Multiplying \eqref{eq:r22} by $\mu_-$ and 
\eqref{eq:r11} by $\mu_+$ and subtracting the 
resulting equations, implies that 
$ad(\mu_-^2 - \mu_+^2)$ is 
meromorphic on $\mathbb S \cap A_s$. Assuming 
$\mu_-^2 - \mu_+^2 = 0$ would imply 
$\mu_{\pm}^4 = 1$ and contradict the assumption 
that $\tr H$ is not independent of $\lambda$. 
Since $\mu_{\pm}$ are 
holomorphic on $\mathbb S \cap A_s$ both $ad$ 
and $bc$ are meromorphic on $\mathbb S \cap A_s$. 
Similarly, evaluating the off--diagonal terms of 
$h\,H \,h^{-1}$ we obtain that both $ab$ 
and $cd$ are meromorphic on $\mathbb S \cap A_s$.

If $d \neq 0$ on some circle 
$C_l$, for $s \leq l \leq r$, then we may write 
$h\,Y^{-1}=\widehat{h}\,\mathrm{diag}[\,1/d,\,d\,]$ 
with
\begin{equation} \label{eq:h-split_case1}
  \widehat{h} := \begin{pmatrix} 
  ad & \tfrac{b}{d} \\ cd & 1 \end{pmatrix} 
  \mbox{ meromorphic on $\mathbb S\cap A_s$.} 
\end{equation}
We define $\mathcal{C} := 
Y^{-1} \,\mathrm{diag}[\,1/d, d\,]\,Y$ and 
$\mathcal{M} := \widehat{h}\, Y$. 
Then $h = \mathcal{M}\, \mathcal{C}$ 
with $\mathcal{C} = \mathcal{C}(\lambda)$ 
defined on $C_l \cap \mathbb S$ and 
$[\,\mathcal{C},\,H  \,] = 0$ there. Moreover, 
$\mathcal{M}$ is meromorphic on 
$\mathbb S \cap A_s$.
 
If $d \equiv 0$ on an arc $C_l \cap \mathbb S$ for 
some $s \leq l \leq r$, then both 
$cd \equiv 0$ on this arc and consequently 
$cd \equiv 0$ in $\mathbb S \cap I_r$, 
and $ad \equiv 0$ on this arc and consequently 
$ad \equiv 0$ in $\mathbb S \cap I_r$.
Thus $d \equiv 0$ on $C_l \cap \mathbb S$ 
implies that both $b,\,c \neq 0$ on 
$\mathbb S \cap I_r$ and we may write 
$h\,Y^{-1} = \widetilde{h}\,
\mathrm{diag}[\,c,\,1/c\,]$ 
with
\begin{equation*} \label{eq:h-split_case2}
  \widetilde{h} = \begin{pmatrix} 
  a/c & bc \\ 1 & cd
  \end{pmatrix} 
  \mbox{ meromorphic on $\mathbb S\cap A_s$.} 
\end{equation*}
Set 
$\mathcal{C} = Y^{-1} \mathrm{diag}
[\, c, 1/c \,] Y$ and   
$\mathcal{M} = \widetilde{h} Y$. Then
$h = \mathcal{M}\, \mathcal{C}$ with 
$[\, \mathcal{C},\,\chi \,] =0$ 
and $\mathcal{M}$ meromorphic on 
$\mathbb S\cap A_s$. \QED


Given an open dense subset $\mathbb S \subset \C^*$, 
obtained by making branch cuts from the points 
$\{ \alpha_j \} \subset \C^*$ to the origin 
respectively $\infty$, depending on whether 
$|\alpha_j|<r$ respectively $|\alpha_j| > r$, let 
$\hat{\mathbb S}$ be the set obtained by making 
branch cuts according to this rule from the points 
$\{ 1/\bar{\alpha}_j \}$. Then the set 
${\mathbb S}^* =  \mathbb S \cap \hat{\mathbb S}$ is 
invariant under $\lambda \mapsto 1/\bar{\lambda}$. 
Note though that if there are at least two distinct 
$\alpha_j \in S^1$, then $\mathbb S^*$ is not 
connected. 

We now improve on these results by exploiting the 
reality condition of loops in $\Lambda_{\C^*} \SU$. 

\mylemma{} \label{th:S}
The eigenvalue functions  $\mu_{\pm}$ of a loop 
$H \in \Lambda_{\C^*}\SU$ 
are of the form $\mu_{\pm} = u \pm iv$, 
where $2u = \tr H$ and 
$v^2 =1 - u^2 = \lambda^N \cdot \hat{v}^2$ with 
$N \in 2\Z$ and $\hat{v}$ holomorphic on 
$\mathbb S$. In particular, 
$S^1 \subset \mathbb S = {\mathbb S}^*$.

\Proof A loop $H \in \Lambda_{\C^*}\SU$ 
is of the form 
\begin{equation*}
  H(\lambda) = \begin{pmatrix} 
  x(\lambda) & y(\lambda) \\
  -\overline{y(1/\bar{\lambda})} & 
  \overline{x(1/\bar{\lambda})} \end{pmatrix}
\end{equation*}
with analytic functions $x,\,y:\C^* \to \C$. Hence
$\left. u \right|_{S^1} \in [-1,\,1]$. 
Consequently, the function $v^2$ is real and 
non-negative on $S^1$ and thus has only even 
roots on $S^1$. If 
$\{ \lambda_1,\,\ldots\,,\lambda_m \}$ are the 
roots of $v^2$ on $S^1$, we may write 
\begin{equation*}
  v^2 = \left( \prod_{j=1}^m 
  (\lambda - \lambda_j)^{2K_j} \right) \, s
\end{equation*}
where $s:\C^* \to \C$ is analytic and 
$\left.s\right|_{S^1} \neq 0$. We write 
$s = c_0\,s_- \, \lambda^N \, s_+$ 
for the scalar Birkhoff decomposition of 
$s$ on $S^1$, with $c_0 \in \C^*$ and $N \in \Z$, 
such that $s_-(\infty)=s_+(0)=1$. 
With these normalizations, the reality condition 
$u^* = u$ implies $s_+^* = s_-$ as well as 
\begin{equation*}
  N = - \sum_{j=1}^m K_j \mbox{ and } 
  \overline{c_0} = c_0 \, \prod_{j=1}^m 
  \lambda_j^{2K_j}.
\end{equation*}
Whence also $s_+^*(\infty)=1$ and 
$v^2 = \lambda^N \hat{v}^2$ with
\begin{equation}
  \hat{v}^2 =  c_0\,s_-\,s_+ \, \prod_{j=1}^m 
  (\lambda - \lambda_j)^{2K_j}.
\end{equation}
Since $s$ is analytic in $\C^*$, both $s_{\pm}$ 
are analytic in $\C^*$, in fact $s_+(\lambda)$ 
is entire, and thus also $\hat{v}^2$ is analytic 
in $\C^*$. There are two cases depending on the 
parity of $N \in \Z$. If $N$ is even, then 
making branch-cuts from the odd ordered 
roots of $s_+$ to the point at infinity and 
branch-cuts from the odd ordered roots of $s_-$ to 
the origin gives an open and dense set 
$\mathbb S \subset \C^*$ on which $\hat{v}^2$ has 
an analytic square root. Since the roots of 
$\hat{v}^2$ on $S^1$ are all even, we have that 
also $v^2$ has an analytic square root on 
$\mathbb S$ and $S^1 \subset \mathbb S$. 
Further, $s_+^* = s_-$ implies 
$\mathbb S^* = \mathbb S$. 

Assume that $N \in \Z$ is odd. Then $\mathbb S$ 
is the set obtained as above but with an 
additional branch cut from the origin to 
the point at infinity. Now $u=\sqrt{\tr^2H-4}$ 
has an absolutely convergent power series expansion 
on $S^1$ and is thus well defined, so there can not 
be a branch cut through $S^1$. Hence $N$ is even.
\QED

\mycorollary{} \label{co:h-split} 
Let $h \in \Lambda_r^+ \SL$ for some $0<r\leq 1$ 
and $H \in \Lambda_{\C^*}\SU$ with non-constant 
trace function and let $\mathbb S$ be as in Lemma 
\ref{th:S}. If $h\, H  \,h^{-1}$ 
is meromorphic on $\mathbb S \cap A_s$ for some 
$0 < s < r$ and in addition we require that 
$(h\, H  \,h^{-1})^* = (h\, H  \,h^{-1})^{-1}$ 
for all $\lambda \in \mathbb S \cap A_s$, then 
\begin{equation}
  h = U\,C 
\end{equation}
with meromorphic loops $U,\,C$ on 
a cut but connected annulus respectively 
a possibly segmented circle $C_l$ for some 
$s<l<r$ with $[\,C,\,H\,]=0$ there, and 
$U \in \Lambda_{r'} \,\SU$ for some 
$r' \in (0,\,1)$.

\Proof
Note that we have $ S^1 \subset \mathbb S = {\mathbb S}^*$ 
under our assumptions.
We apply Theorem \ref{th:h-split} and decompose  
$h = \mathcal{M} \,\mathcal{C}$ on $C_l$ for some 
$s<l<r$ with $\mathcal{M}$ meromorphic on 
$\mathbb S \cap A_s$ and 
$[\, \mathcal{C} ,\,H  ]=0$.  
Consequently, we have $h\,H \,h^{-1} = 
\mathcal{M}\,H\,\mathcal{M}^{-1}$ on 
$\mathbb S \cap A_s$ and 
$(h\, H  \,h^{-1})^* = 
  (h\, H  \,h^{-1})^{-1}$ 
is equivalent to 
\begin{equation}
  [\,\mathcal{M}^* \,\mathcal{M} ,\,
  H  \,]=0
\end{equation}
on $\mathbb S\cap A_s$. Next, we seek a loop $L$, 
meromorphic on a subset of $\mathbb S \cap A_s$, 
of the form  $L = f \, \Id + g \, \mathcal{M}^*\mathcal{M}$ 
with meromorphic functions $f,\,g$ on 
$\mathbb S \cap A_s$, such that 
\begin{equation} \label{eq:Luni}
  (\mathcal{M} L^{-1})^* = 
  (\mathcal{M} L^{-1})^{-1}.
\end{equation}
Then $h = \mathcal{M} L^{-1} L \mathcal{C}$ 
is the desired factorization. 
Equation \eqref{eq:Luni} is equivalent to 
$\mathcal{M}^* \mathcal{M} = L^* L$ and by setting 
$P := \mathcal{M}^* \mathcal{M}$, yields 
\begin{equation} \label{eq:P} \begin{split}
  P &= ff^* \Id + (f^* g + f\,g^*)P + g\,g^* P^2 \\ 
  &= (ff^* - gg^*)\Id + (f^*g + fg^* + gg^* \tr(P))P
  \end{split}
\end{equation}
since $P^2 = \tr(P)P - \Id$ by Cayley-Hamilton. 
The Ansatz $f=f^*=g$ reduces \eqref{eq:P} to 
\begin{equation}\label{eq:g^2}
   g^2 = \frac{1}{2 + \tr P}.
\end{equation}
If we denote the entries of $\mathcal{M}$ by 
$m_{ij}$, then the function
\begin{equation}
  \tr P(\lambda) = \sum_{i,j=1}^2 m_{ij}(\lambda)
  {m_{ij}^*(\lambda)}
\end{equation}
is meromorphic on $\mathbb S \cap A_s$ and 
is real and positive on $S^1\subset \mathbb S \cap A_s$ .
With the help of a scalar Birkhoff decomposition of $g^2$ 
it is straightforward to see that
there exists a well defined square root of \eqref{eq:g^2} 
on $S^1$, which extends meromorphically to the set 
$\mathbb S_P \subset \mathbb S$, obtained by making branch-cuts 
from the odd ordered roots 
of $2 + \tr P$ in the usual fashion. 
Note that $S^1 \subset \mathbb S_P$ and 
$\mathbb S_P^* = \mathbb S_P$. 
Now that $g = (2 + \tr P)^{-1/2}$ is 
meromorphic on $\mathbb S_P \cap A_s$, we can set 
$U:=\mathcal{M}L^{-1}$, which is meromorphic on 
$\mathbb S_P \cap A_s$ with $U^* = U^{-1}$, 
and $C : = L\,\mathcal{C}$, which is defined on 
$C_l \cap \mathbb S_P$ with $[\,C,\,H\,]=0$. 
It remains to show that there exists 
an $r' \in (0,\, 1)$ such that 
$U \in \Lambda_{r'} \,\SU$. 

We choose $0<r'<1$ such that both of the 
following conditions hold:

(i) The only poles of $\mathcal{M}$ in 
    $A_{r'} \cap \mathbb S$ are on $S^1$, 

(ii) $A_{r'} \cap \mathbb S_P = A_{r'}$.

Denote the pole set of $\mathcal{M}$ in 
$A_{r'}$ by 
$\mathcal{P} = \{ p_1,\dots,p_K \} \subset S^1$. 
Let $V_j \subset A_{r'}$ be an open 
neighborhood of $p_j$ such that 
$V_j \cap \mathcal{P} = \{ p_j \}$. 
Let $V^*_j = V_j \setminus \{ p_j \}$.
Since $U$ is unitary on 
$\mathbb S_P \cap U^*_j \cap S^1$,  
its entries are bounded. 
Hence $U$ extends analytically to $p_j$ 
for all $j=1,\dots,K$. Consequently, 
$U$ is analytic in $A_{r'}$ with $U^* = U^{-1}$. 
Thus $U \in \Lambda_{r'}\,\SU$, concluding this proof. 
\QED

In the following sections we will apply the 
results of the present section to the twisted 
loop groups.


\newsection{} 
Let $F \in \FID$ and $h \in \LGP$ and let
$\hat F = h \# F$. Then $\hat F \in \FID$. 
Let $\tau \in \Delta$ and $\chi(\tau)$, 
$\hat\chi(\tau)$ be monodromies of $F$ 
respectively $\hat F$. Writing $hF = \hat FB$ 
for an Iwasawa decomposition and 
$\tau^*F = \chi(\tau)\,F\,k(\tau)$ and 
$\tau^*\hat F=\hat\chi(\tau)\,\hat F\,\hat k(\tau)$, 
we obtain 
\begin{equation} \label{eq:mono_dress1}
  \hat\chi(\tau)^{-1} h\,
  \chi(\tau) \, h^{-1} \hat F =
  \hat F\,\hat k(\tau)\,(\tau^* B)\,
  k(\tau)^{-1}B^{-1}.
\end{equation}
Since $\hat k(\tau)\,(\tau^* B)\,k(\tau)^{-1}
B^{-1} \in \gauge$, 
this proves that 
\begin{equation} \label{eq:mono_dress2} 
  \hat \chi (\tau)^{-1}h\,\chi(\tau)\,h^{-1} 
  \in \mathrm{Iso}_r(\hat F).
\end{equation}
Since $F(\tilde{z}_0)= \hat F(\tilde{z}_0)=\Id$, 
evaluation of $hF = \hat FB$
at the base-point $\tilde{z}_0$ yields 
$h=B(\tilde{z}_0)$. Evaluating equation 
\eqref{eq:mono_dress1} at $\tilde{z}_0$ 
and solving for $\hat \chi(\tau)$ gives  
\begin{equation} \label{eq:monodromy_dressed1}  
  \begin{split}
  \hat \chi (\tau) &= h\,\chi(\tau)\, 
  k(\tau,\tilde{z}_0)\,
  B\left(\tau(\tilde{z}_0)\right)^{-1}
  \hat k(\tau,\tilde{z}_0)^{-1} \\
  &= h\, \chi(\tau) \, h^{-1} \,L,\,
  L \in \mathrm{Iso}_r(\hat F) \\
  &= h\,\chi(\tau)\,b\, h^{-1},\,
  b := h^{-1} L \, h \in \mathrm{Iso}_r(F). 
  \end{split}
\end{equation}
This leads us to an application of Theorem 
\ref{th:h-split} and shows that if a 
dressing matrix preserves the topology, then it 
factorizes as in \eqref{eq:h-split}.  

We now denote by $\mathbb S = {\mathbb S}^*$ 
the set obtained from the eigenvalues of 
a monodromy $\chi$.

\mytheorem{} \label{th:h-split1} 
Let $\tau \in \Delta$ and $F \in \FID$ 
with monodromy $\chi = \chi(\tau)$. We 
assume that the trace of $\chi$ is not 
constant. Let $h \in \LGP$ such that 
$h\,\chi\,h^{-1}$ is meromorphic on 
$\mathbb S \cap A_s$ for some $0<s<r$, 
where $\mathbb S$ is as in Lemma \ref{th:S}. 
If $\hat F = h\#F$ has a monodromy 
$\hat{\chi}$ with respect to 
$\tau$, and if $\hat{\mathbb S}$ is defined 
for $\hat{\chi}$ according to Lemma \ref{th:S},
then $h$ can be factored, on a possibly 
segmented circle 
$C_l \cap \mathbb S \cap \hat{\mathbb S}$ for 
$l \in [s,\,r]$, into 
$h = \mathcal{M} \,\mathcal{C}$ 
where $\mathcal{C} = \mathcal{C}(\tau,\,\lambda)$ 
is twisted $\SL$--valued and defined on 
$C_l \cap \mathbb S \cap \hat{\mathbb S}$ and 
$[\, \mathcal{C},\, \chi  \,] = 0$ there, and 
$\mathcal{M} = \mathcal{M}(\tau,\,\lambda)$ 
is twisted $\SU$--valued and meromorphic 
on the connected open set 
$ S^1 \subset \mathbb S \cap 
\hat{\mathbb S} \cap A_s$.

\Proof By the above \eqref{eq:monodromy_dressed1}, 
for the dressed monodromy we have  
$\hat{\chi} = h\,\chi\,h^{-1}\,\hat{L}$ with 
$\hat{L} \in \mathrm{Iso}_r(\hat F)$. Clearly, 
$\hat{L}$ is defined and meromorphic on 
$\mathbb S \cap A_s$. Moreover, we can write 
$\hat{\chi} = h\,\chi\, L h^{-1}\,$, where 
$L \in \mathrm{Iso}_r(F)$. 
From this we conclude that the eigenvalue 
functions of $\chi\,L$ are the same as those 
of $\hat{\chi}$. At this point we undress 
every occurring matrix and continue, 
until the very end of the proof, to work 
in the undressed setting. We note that 
the domains of definition of the dressed and 
the undressed matrices are in a natural 
correspondence and that unitary matrices 
are mapped to unitary matrices, 
while the "positive matrices" almost correspond. 
For simplicity of notation we will use the 
same notation as in the twisted situation.
First we note that we can apply Theorem 
\ref{th:h-split} to $h\,\chi\,h^{-1}$.
Hence  $h = \mathcal{M} \,\mathcal{C}$, 
where $\mathcal{M}$ has values in $\SL$ 
and is defined on $\mathbb S$.
As a consequence, $\hat{\chi} = 
\mathcal{M} \chi {\mathcal{M}}^{-1} \hat{L} = 
\mathcal{M} \chi \,L'{\mathcal{M}}^{-1}$. 
From this expression we infer that $L'$ is 
defined and meromorphic on $\mathbb S$.
Now we can apply almost verbatim the proof of 
\ref{co:h-split}, if we replace everywhere
$\mathbb S$ by $\mathbb S \cap \hat{ \mathbb S} = 
\tilde{\mathbb S}$ and obtain that we can assume 
that $\mathcal{M}$ is defined and holomorphic 
on some open neighborhood of $S^1$ and unitary 
on $S^1$. Finally, twisting every occurring matrix 
again we obtain the desired result. \QED

When the associated family of $F \in \FID$ 
possesses umbilic points, then, since dressing 
preserves umbilic points \cite{Wu:dre}, 
$\rm{Iso}_r(\hat F) = \left\{ \pm \Id \right\}$ 
for $\hat F = h \# F$ and \eqref{eq:mono_dress2} 
implies 
\begin{equation} \label{eq:monodromy_dressed2}
  \hat \chi (\tau) = \pm h\,\chi(\tau)\,h^{-1} 
  \in \Lhol.
\end{equation}
This allows us to restate Corollary 
\ref{co:h-split} for the trivial isotropy case: 


\mycorollary{} \label{co:h-split1} 
Let $F \in \FID$ with 
$\mathrm{Iso}(F) = \{ \pm \Id \}$. Let 
$\chi = \chi(\tau) \in \Lhol$ be a
monodromy of $F$ with non-constant trace. 
Let $h \in \LGP$ for some $r \in (0,\,1]$. 
If $h\#F$ also has a 
monodromy with respect to $\tau$, then $h$ can be 
factored into $h = U \,C$ 
with $[\, C,\,\chi  \,] = 0$  and 
$U = U(\tau,\,\lambda)$ the meromorphic 
extension of an element of 
$\Lambda_{r^\prime} \SU_\sigma$ 
for some $r^\prime \in [r,\,1]$. \QED

We were not able to prove the corresponding result for a 
general $F \in \FID$ without the isotropy 
condition. It turns out, that we do 
have an analogous result for the dressing orbit of 
the vacuum, to which we now turn our attention.


\section{Dressing the vacuum}
\newsection{\bf Definitions.} 
The standard round cylinder (=''the vacuum'') 
has the conformal structure of the punctured 
plane $\C^*$. If we identify 
$\C^* \cong \C / q\, \Z$ for some 
$q \in \C^*$, then 
\begin{equation}
  \exp:\C \to \C^*,\,w \mapsto z = \exp(q\,w)
\end{equation}
is the universal covering map. The group of deck 
transformations $\Delta \cong \Z$ is generated by 
the translation
\begin{equation}  \label{eq:translation}
  \tau_q : w \mapsto w + q.
\end{equation}
An extended unitary frame of the associated 
family of the vacuum is given by
\begin{equation} \label{eq:F_c}
  F_c(w,\,\lambda) = 
  \exp((w \lambda^{-1} 
  - \overline{w} \lambda)\, A)
\end{equation}
where
\begin{equation}
  A = \begin{pmatrix} 
  0 & 1 \\ 1 & 0 \end{pmatrix}.
\end{equation}
A monodromy with respect to $\tau_q$ of 
$F_c$ is given by 
\begin{equation} \label{eq:cylinder_monodromy}
  \chi_c (\tau_q) = 
  \exp((q \lambda^{-1} - \bar{q} 
  \lambda )\,A).
\end{equation}
The map $F_c:\C \to \Lhol$ is obtained in the 
DPW framework from the triple 
$((\lambda^{-1} + \lambda )A\,dw,\,\Id,\,0 )$. 
Note that we have chosen the base-point 
$w_0 = 0$. The corresponding holomorphic frame 
is 
\begin{equation} \label{eq:Phi_c}
  \Phi_c(w,\,\lambda) = 
  \exp((\lambda^{-1} + \lambda)\,w\, A)
\end{equation}
with monodromy with respect to $\tau_q$ given by 
\begin{equation} \label{eq:rho_c}
  \varrho_c(\tau_q) = 
  \exp((\lambda^{-1} + \lambda)\,q\, A).
\end{equation}
Note that if $q \in i\R$, then 
$\varrho_c(\tau_q) \in \Lhol$ and 
$\varrho_c(\tau_q) = \chi_c(\tau_q)$. 

%
\newsection{\bf Commuting flows.} 
Denote the abelian sub-algebra of elements of 
$\Lambda_{\mbox{\tiny{$\C^*$}}}\sl_\sigma$ 
that commute with $A$ and have a pole at 
$\lambda = 0$ by 
\begin{equation*}
  \mathcal{Z} = \left\{ \phi \,A: 
  \phi (\lambda) = \Sigma
  \phi_j \lambda^j \mbox{ odd },\, j\geq -K,\, 
  K = 2g +1 \in \N \right\}. 
\end{equation*}
Let $h \in \LGP,\,\eta \in \mathcal{Z}$ and 
define $h \sharp \eta$ via the $r$--Iwasawa 
decomposition 
\begin{equation}
  h \exp(\eta) = U \, h \sharp \eta.
\end{equation}
Further, for $F = h \# F_c$ define 
$F \sharp \eta := (h \sharp \eta) \# F_c$. 
By \cite{BurP:dre}, Proposition 4.1, this is 
an action on the dressing orbit of $F_c$: 
If $\eta_1,\,\eta_2 \in \mathcal{Z}$, then 
$F \sharp (\eta_1 + \eta_2) = 
(F \sharp \eta_2)\sharp\eta_1  = (F \sharp \eta_1)\sharp\eta_2$ and 
if $\eta_1 - \eta_2 \in \Lambda^+_r \sl$, then 
$F \sharp \eta_1 = F \sharp \eta_2$. 
By construction $\zeta = \phi A \in \mathcal{Z}$ 
is analytic for $\lambda \in \C^*$ with a pole at 
$\lambda = 0$. Expanding $\phi(\lambda) = 
\sum_{j\geq -K} \phi_j \lambda^j$, observe that 
$\exp(\zeta) \# F_c = \exp(\eta) \# F_c$ 
where $\eta = \phi^\prime A$ and 
$\phi^\prime = \sum_{j=-K}^0\phi_j \lambda^j$ is 
analytic on $\CP \setminus \{ 0 \}$. Since  
$\eta^* = \phi^* A$ is analytic 
in $\C$, $e^{\eta - \eta^*} \in \Lhol$ and 
$e^{\zeta}\# F_c = e^{\eta-\eta^*}F_c$. 
Hence, when dressing with $e^\zeta$, we restrict 
without loss of generality to 
$\zeta \in \mathcal{Z}$ of the form 
\begin{equation} \label{eq:zeta}
  \zeta = \sum_{j=1}^K (\phi_j \lambda^{-j} - 
  \bar{\phi}_j \lambda^j ) A
\end{equation}
where $\phi_j \in \C$ and summation is over odd 
indices $j=1,\,3,\ldots,K = 2g+1 \in \N$. Then
\begin{equation}
  \exp (\zeta) \, \# \, 
  F_c(w,\,\lambda) = 
  \exp \bigl(\sum_{j=3}^{K} (\phi_j \lambda^{-j} - 
  \bar{\phi}_j \lambda^j )A \bigr) \, 
  F_c(w+\phi_1,\,\lambda).
\end{equation}
%

For $\zeta \in \mathcal{Z}$ of the form 
\eqref{eq:zeta} we have $e^{t\zeta}\in \Lhol$ 
for all $t \in \R$. 
Hence, for any extended unitary frame 
$F \in \mathcal{F}(\MT)$, we have 
$e^{t\zeta}\#F = e^{t\zeta}F$ which on the level 
of the immersions has the effect 
\begin{equation*}
  f \longmapsto \exp(t\,\zeta)\,f\,\exp(-t\,\zeta) - 
  \tfrac{1}{2H} i\lambda t 
  \tfrac{\partial \zeta}{\partial \lambda}.
\end{equation*}
Each $\zeta \in \mathcal{Z}$ defines a flow on the 
dressing orbit of the standard round cylinder by 
\begin{equation}
  F_t := F \sharp t \zeta = (h\sharp t\zeta)\#F_c 
  = U(\lambda,\,t)^{-1} h\#e^{t\zeta}F_c
\end{equation}
where we have defined $U(\lambda,\,t)$ via the 
$r$--Iwasawa decomposition 
$he^{t\zeta}=U(\lambda,\,t)(h\sharp t\zeta)$ and 
$F = h \# F_c$ for some $h \in \LGP$. 
The flow $F_t$ is called 
\textit{trivial} if and only 
if $F_t = V(t) F V(t)^{-1}$ for a 
$\lambda$--independent map $V:\R \to \SU$ 
and $F$ is of \textit{finite type} if and 
only if the subspace 
$\mathcal{Z}^\prime \subset \mathcal{Z}$ of 
trivial flows has finite co-dimension in 
$\mathcal{Z}$.


\newsection{\bf Monodromy.} 
Let $h \in \LGP$ and $hF_c = FB$ the 
$r$--Iwasawa decomposition of $hF_c$ such that 
at the  base-point $w_0=0$ we have 
$F_c(0,\,\lambda) = F(0,\,\lambda) = \Id$ and 
$h(\lambda) = B(0,\,\lambda)$. 
Let us assume that $F$ has symmetry 
$(\tau_q,\,\chi(\tau_q))$. 
From \eqref{eq:monodromy_dressed1} we obtain
\begin{equation}\label{eq:dress_hol} 
  \chi (\tau_q) = h\, \chi _c(\tau_q)\, B(q)^{-1}.
\end{equation}
For completeness of exposition we prove the 
following result, first obtained in sections 
3.1--3.4 of \cite{DorH:per}. 

\mylemma{}\label{th:f_q}
Let $h \in \LGP$ and $F_c$ 
be as in \eqref{eq:F_c}. 
Assume that $h \# F_c$ has monodromy 
$\chi(\tau_q) \in \Lhol$. 
Then there exists an odd function 
$f_q :C_r \to \C$ with a holomorphic extension 
to $I_r$ such that 
$B(q)^{-1} = \exp(f_q A) \, h^{-1}$. 
Setting 
\begin{equation} \label{eq:p_q}
  p_q(\lambda) = q\lambda^{-1}-\bar{q}\lambda 
  + f_q(\lambda),
\end{equation}
we may write the monodromy of $h \# F_c$ as 
\begin{equation} \label{eq:dress_monodromy}
  \chi (\tau_q) = h\, \exp(p_q A)\, h^{-1}.
\end{equation}

\Proof 
Let $\chi_c(\tau_q)$ be the monodromy of 
$F_c$ as in \eqref{eq:cylinder_monodromy}.
Then
\begin{equation}
  \chi_c(\tau_q) F_c = \tau_q^* F_c = 
  h^{-1} \chi (\tau_q) h F_c B^{-1} \tau_q^* B
\end{equation}
in combination with \eqref{eq:dress_hol} gives 
$B(q)^{-1}B(0)  F_c = F_c \, \tau_q^* B^{-1} B$. 

Hence $[B(q)^{-1}B(0), A] = 0$ and for the matrix 
\begin{equation}\label{eq:T}
  T = \frac{1}{\sqrt{2}}  \begin{pmatrix} 
  1 & 1 \\ -1 & 1 \end{pmatrix}
\end{equation}
we have $TAT^{-1} = \sigma_3$ and 
$[TB(q)^{-1}B(0)T^{-1}, \sigma_3] =0$. 
Thus 
$TB(q)^{-1}B(0)T^{-1} = 
\mathrm{diag}[s_q,\, s_q^{-1}]$ 
for some holomorphic function $s_q:I_r \to \C^*$. 
Defining $s_q = \exp(f_q)$ yields 
$B(q)^{-1}B(0) = \exp(f_q A)$. 
Since $B(0) = h$, this proves the claim.
\QED

The point of Lemma \ref{th:f_q} is that 
we can now write a monodromy of an extended frame 
in the dressing orbit of the vacuum as
\begin{equation}\label{eq:sinhcosh}
  \cosh(p_q) \,\Id + \sinh(p_q) \,h\,A\,h^{-1}
\end{equation}  
and much of the analysis reduces to the study of 
the two scalar functions
\begin{equation} \label{eq:alphabeta}
  \alpha=\cosh(p_q)\mbox{ and }\beta = \sinh(p_q).
\end{equation}
In analogy to Corollary \ref{co:h-split1} 
we have the following result for the monodromy 
representation of surfaces in the dressing orbit 
of the cylinder.

\mytheorem{}\label{th:h-split_cyl}
Let $h \in \LGP$ and $F_c$ 
be as in \eqref{eq:F_c}. 
Assume that $h \# F_c$ has symmetry 
$(\tau_q,\,\chi(\tau_q))$. 
Then $h = \mathcal{M}(\tau_q)\,\mathcal{C}(\tau_q)$ 
with $\mathcal{M}(\tau_q) \in \Lambda_{r^\prime} 
\SU_\sigma$ for suitable $r \leq r^\prime < 1$ and 
$[\,\mathcal{C}(\tau_q),\,A\,]=0$.

\Proof For $T$ defined in \eqref{eq:T} and by
Lemma \ref{th:f_q}, we may write
\begin{equation} \label{eq:chi}
  \chi(\tau_q) = h \,T^{-1} \,
  \mathrm{diag}[\,\exp(-p_q),\, 
  \exp(p_q)\,]\, T\, h^{-1}.
\end{equation} 
for $0< \mid \lambda \mid < r$. 
Let $h T^{-1} = \bigl( \begin{smallmatrix} 
a & b \\ c & d \end{smallmatrix} \bigr)$.
The diagonal entries of $\chi(\tau_q)$ 
are analytic functions in $\lambda \in \C^*$ 
and given by
\begin{equation} \begin{split} 
  a \, d \, \exp(-p_q) - b\,c\,\exp(p_q),\\
  a \, d \, \exp(p_q) - b\,c\,\exp(-p_q). \end{split}
\end{equation}
Adding these implies that $\alpha = \cosh(p_q)$ 
is an analytic function on $\C^*$. 
More precisely, $\alpha$ has an analytic 
extension from $ 0 <\mid \lambda \mid <r$ to $\C^*$.
Solving for $\exp(p_q)$ we obtain
%
%
 $\exp(p_q) = \alpha + \sqrt{\alpha^2-1}$. 
This is analytic on 
the set $\mathbb S \subset \C^*$ obtained by 
making branchcuts from the odd--ordered roots of 
$\alpha^2-1$ in the usual fashion. 
Note that the equation \eqref{eq:chi} above 
implies that $\mu_{\pm} = \exp (\pm p_q)$ 
are the eigenvalues of $\chi (\tau_q)$. 
Hence, by Lemma \ref{th:S} we know that 
$\mathbb S$ is open, connected, and dense 
in $\C^*$. Moreover we have $S^1 \subset \mathbb S$ 
and $\mathbb S = {\mathbb S }^*$.

Consequently, $\exp(-p_q)$ is holomorphic on 
$\mathbb S$. Multiplying the $(1,1)$ entry of 
$\chi(\tau_q)$ by $\exp(p_q)$ and the $(2,2)$ 
entry of $\chi(\tau_q)$ by $\exp(-p_q)$ 
and subtracting the resulting equations implies 
that $bc \sinh (p_q)$ is holomorphic on $\mathbb S$. Note that $\beta = \sinh (p_q) \not\equiv 0$. Hence
$bc$ is meromorphic on $\mathbb S$. 
Similarly one shows that $ad$ is meromorphic on 
$\mathbb S$. 

The off--diagonal terms of $\chi(\tau_q)$ 
are $2 ab \beta$ and $-2 cd \beta$. 
Hence both $ab$ and $cd$ are meromorphic on 
$\mathbb S$. 
If $d \neq 0$ on the circle $C_r$, then as in 
\eqref{eq:h-split_case1}, we may write 
$hT^{-1} = \hat{h}D_1$ with 
$D_1= \mathrm{diag}[1/d,\,d]$ 
and $\hat{h}$ meromorphic on $\mathbb S$. 

Defining $\mathcal{C}(\tau_q) = T^{-1}D_1 T$, 
$\mathcal{M}(\tau_q) = \hat{h} T$ yields 
$h = \mathcal{M}(\tau_q) \mathcal{C}(\tau_q)$ 
with $[\, \mathcal{C}(\tau_q),\,A \,] =0$ and 
$\mathcal{M}(\tau_q)$ meromorphic on $\mathbb S$.
Arguing as in the proof of Theorem \ref{th:h-split}, 
$d \equiv 0$ on an arc $C_s \cap \mathbb S$ 
implies that both $b,\,c \neq 0$ on $\mathbb S$ 
and we may write $hT^{-1} = \widetilde{h}D_2$ with 
$D_2= \mathrm{diag}[c,\,1/c]$ and 
$\widetilde{h}$ meromorphic on $\mathbb S$. 
Defining $\mathcal{C}(\tau_q) = T^{-1}D_2 T$, 
$\mathcal{M}(\tau_q) = \widetilde{h} T$ yields 
$h = \mathcal{M}(\tau_q) \mathcal{C}(\tau_q)$ 
with $[\,\mathcal{C}(\tau_q),\,A\,] =0$ and 
$\mathcal{M}(\tau_q)$ meromorphic on $\mathbb S$ 
in this case. In either case, on $C_r$ we may write
\begin{equation} \label{eq:U} 
  \chi(\tau_q) = \mathcal{M}(\tau_q)\, 
  \exp(p_q\,A)\,\mathcal{M}(\tau_q)^{-1}.
\end{equation}
Introducing an additional cut from $0$ to $-\infty$ 
we obtain a simply connected domain 
$\hat{\mathbb S} = \hat{\mathbb S}^*$, 
on which we can take the logarithm of 
the holomorphic function $\exp(p_q)$, 
thus defining $p_q$ on $\hat{\mathbb S}$. 
Then on $\hat{\mathbb S}$ we combine 
$\chi(\tau_q)^* = \chi(\tau_q)^{-1}$ 
with \eqref{eq:U} to obtain
\begin{equation} \label{eq:mstarm}
  \mathcal{M}^*(\tau_q) \,\mathcal{M}(\tau_q) 
  \,\exp(-p_q\,A) = \exp(p^*_q\,A)\, 
  \mathcal{M}^*(\tau_q) \,\mathcal{M}(\tau_q).
\end{equation}
For $T$ defined in \eqref{eq:T}, recall that 
$T\,A\,T^{-1} = \mathrm{diag}[1,\,-1]$. 
Since $T \in \SU$, the matrix  
$X:= T\,\mathcal{M}^*(\tau_q) \,
\mathcal{M}(\tau_q)\,T^{-1}$ is hermitian and 
of the form $X = W\,W^*$. In particular, if we write 
$X = \bigl( \begin{smallmatrix} x & y \\ 
\bar{y} & z \end{smallmatrix} \bigr)$, 
then $x \neq 0$.
We can rewrite \eqref{eq:mstarm} as 
$X\,T\,\exp(-p_q\,A)\,T^{-1} = 
T\,\exp(p^*_q\,A)\,T^{-1}\,X$ and obtain
\begin{equation}
  \begin{pmatrix} x & y \\ \bar{y} & z \end{pmatrix}
  \begin{pmatrix} \exp(-p_q) & 0 \\ 0 &  \exp(p_q) 
  \end{pmatrix} = \begin{pmatrix} \exp(p^*_q) & 0 \\ 
  0 & \exp(-p^*_q) \end{pmatrix} 
  \begin{pmatrix} x & y \\ \bar{y} & z \end{pmatrix}.
\end{equation}
Since $x \neq 0$, then $\exp(-p_q) = \exp(p^*_q)$ 
on $\hat{\mathbb S}$ and we conclude that 
$\cosh(p^*_q) = \cosh(p_q)$ and 
$\sinh(p^*_q) = - \sinh(p_q)$. Consequently 
$\exp(p^*_q\,A) = \exp(-p_q\,A)$ 
on $\hat{\mathbb S}$ and equation 
\eqref{eq:mstarm} reads 
\begin{equation}
  [\,\mathcal{M}^*(\tau_q) \mathcal{M}(\tau_q),\,
  \exp(-p_q\,A)\,] = 0.
\end{equation}
From the proof of Corollary \ref{th:S} we 
now obtain some $L(\tau_q)$ 
and $r' \in [r,1)$ such that 
$\mathcal{M}(\tau_q)L(\tau_q) \in 
\Lambda_{r'}\SU_\sigma$ and 
$[\,L^{-1}(\tau_q),\,\mathcal{C}(\tau_q)\,]=0$, and 
concludes the proof. \QED
%


Before going on to evaluate the result above 
to examples we would like to relate the work 
above to \cite{DorH:per}. Of particular 
interest to us is equation (3.6.7) there:
\begin{equation} \label{eq:x}
  h_+ = \frac{i}{2 \sqrt{\tilde{b}}}
  \begin{pmatrix}
  (x-x^{-1}) \tilde{b} & (x+x^{-1}) \tilde{b} \\
  (x+x^{-1}) - (x-x^{-1}) \tilde{a} & 
  (x-x^{-1}) - (x+x^{-1}) \tilde{a} \end{pmatrix}
\end{equation}
where $\tilde{a}$ and  $\tilde{b}$ 
denote the entries of the matrix $hAh^{-1}$ 
in the $(11)$--position and the $(12)$--position 
respectively.
This is the general form of a dressing matrix 
having on $\chi$ the same effect as $h$.
The question addressed in the Theorem above 
is thus equivalent to the question whether one 
can choose $x$ so that $h_+$ is unitary on $S^1$.

\mylemma{} The matrix $h_+$ is unitary on 
$S^1$ if and only if
\begin{equation}
  |\,x\,|^2 = \frac{1 + \tilde{a} +  
  \sqrt{|\tilde{b}|^2}}{1-\tilde{a} + 
  \sqrt{|\tilde{b}|^2}}. 
\end{equation}
This equation can be solved with some $x$ 
holomorphic at $\lambda =0$ and well defined 
on some sufficiently well cut complex plane.

\Proof
The unitarity condition for $h_+$ yields 
two equations. However, a straightforward 
computation shows that
these two equations are equivalent. 
It thus remains to consider
\begin{equation}
  (x^* - (x^*)^{-1}) \sqrt{\tilde{b}^*} = 
  \frac{1}{\sqrt{\tilde{b}}} \lbrack (x-x^{-1}) - 
  (x+x^{-1}) \tilde{a} \rbrack.
\end{equation}
Splitting this equation into real and imaginary 
part yields two equations. 
However, again, a straightforward computation 
shows that these two equations are equivalent.
It thus suffices to consider
\begin{equation}
  (Re(x)^* - (Re(x)^*)^{-1}) \sqrt{|\tilde{b}|^2} =  
  (Re(x)-Re(x)^{-1}) - (Re(x)+ Re(x)^{-1}) \tilde{a}.
\end{equation} 
Note that we have used here that $\tilde{a}$ 
is real on $S^1$ by \cite{DorH:per}, Theorem 3.7.
Using $x^{-1} = x/{\mid x \mid}^2$ this 
equation rewrites directly into the claim.
\QED

\mycorollary{}
$h_+$ can be chosen diagonal if and only 
if $\tilde{a} =0$. In this case the 
diagonal entries of $h_+$ are unitary on $S^1$.

\Proof
Using \eqref{eq:x} for general $x$ we see that 
$h_+$ is diagonal only if 
$x + x^{-1}$ vanishes identically, 
since $\tilde{b}$ cannot vanish by b) and c) of 
\cite{DorH:per}, Theorem 3.7.
Inserting this into the second off-diagonal 
expression for $h_+$ we obtain $\tilde{a} =0$. 
The converse is obvious choosing $x = i$.
Assume now that $h_+$ has been chosen as a 
diagonal matrix. Since  $\tilde{a}$ vanishes 
in this case as just shown, c) and e) of 
\cite{DorH:per} imply that   
$\tilde{b}$ is unitary on $S^1$.
Finally, in view of \eqref{eq:x}, we obtain
\begin{equation}
  h_+ = \frac{i}{2 \sqrt{\tilde{b}}}(x - x^{-1})
        \mathrm{diag}[b,\,1].
\end{equation}
Moreover, $x = -x^{-1}$, whence $x = \pm i$. 
The condition $\det h_+ = 1$ is now satisfied.
Thus necessarily 
$h_+ = \pm \mathrm{diag}[\sqrt{\tilde{b}},\,
1/\sqrt{\tilde{b}}]$ 
and the diagonal entries of $h_+$ are unitary 
on $S^1$. \QED
 
An immediate consequence of the above discussion leads to 
following observation.

\mycorollary{}
Let $h \in \LGP$ such that 
$h\, \exp(p_q A)\, h^{-1} \in \Lhol$ for some 
translation $\tau_q$. If $h$ is diagonal, 
then $h^4$ is rational on $\CP$ and
$h^2$ is analytic in 
$\CP \setminus \mathcal{H}$, where $\mathcal{H}$ 
is the set obtained by making appropriate 
branchcuts from the points lying in the set
\begin{equation} \label{eq:branchset}
  \mathcal{H}^\prime = 
  \{ \lambda \in \CP : |\lambda | > r \mbox{ and } 
  p_q(\lambda) \in \pi i \Z \} \cup \{ 
  \mbox{ singularities of } p_q(\lambda) \}
\end{equation}
and $h$ is unitary on $S^1$. \QED

\myremark{}
Writing $h = \mathrm{diag}[a,\,1/a]$, we have
that $\chi(\tau_q)= \alpha\,\Id + 
\beta\,\mathrm{off}[a^2,\,a^{-2}]$ 
is analytic on $\C^*$. 
Hence $a^2$ and $a^{-2}$ are analytic in 
$\mathcal{H}$ and even in $\lambda$.
Following the procedure for splitting 
a matrix we first untwist $h$. 
The above procedure for 
constructing the factorisation 
$h = \mathcal{MC}$ gives
\begin{equation} \label{eq:diagonal_h-split}
  \mathcal{M} = \frac{1}{2} \begin{pmatrix} 
  1+a^2 & 1-a^2 \\
  -1 + a^{-2}  & 1 + a^{-2} \end{pmatrix}, \,
  \mathcal{C} = \frac{1}{2} \begin{pmatrix} 
  a+1/a & a-1/a \\ a-1/a & a+1/a \end{pmatrix},
\end{equation}
where a really denotes the untwisted function.
Now we need to twist again. This way we obtain 
the original function $a$ and in addition 
the off-diagonal terms are multiplied by 
$\lambda$ (in the $(12-)$position) and 
by $\lambda^{-1}$ in the $(21)-$position 
respectively. Clearly this leads to a matrix 
$\mathcal{M}$ which is twisted as required. 


\newsection{\bf Delaunay Surfaces.} 
\label{sec:Delaunayexample} 
As an example, we discuss how to dress the 
vacuum into a Delaunay surface and compute 
the splitting of this dressing matrix according 
to Theorem \ref{th:h-split_cyl}. It is proven in 
\cite{Kil:del} that an associated 
family of the Delaunay surface with neck radius 
$\omega$ is generated by the triple 
$(D\, dz/z,\,\Id,\,0 )$, where
\begin{equation} \label{eq:delaunay}
  D = \begin{pmatrix} 0 & \alpha \lambda^{-1} + 
  \beta \lambda \\ \beta \lambda^{-1} + 
  \alpha \lambda & 0 \end{pmatrix},
\end{equation}
and $\alpha,\, \beta \in \R$ with 
$\alpha + \beta = 1/2$ and
$\omega = \tfrac{1}{2H}
(1 - \sqrt{1 - 16 \alpha \beta})$. 
The monodromy of the unitary frame with respect to 
the translation $\tau_q$ for $q = 2 \pi i$ is 
given by $\exp ( 2 \pi i D )$. 
Therefore, to dress the 
standard cylinder into a Delaunay surface, 
it suffices to determine a diagonal 
$h = \mathrm{diag}[a,\, 1/a] \in \LGP$
for some suitable $r>0$ such that 
$p_q h A h^{-1} = D$, which is equivalent to 
the two equations
\begin{equation*}
  a^{\pm 2}(q\lambda^{-1}- \bar{q}\lambda + f_q )= 
  2 \pi i  (\alpha\lambda^{\mp 1} + 
  \beta \lambda^{\pm 1} ).
\end{equation*}
Solving both equations for $a^2$ and equating gives
\begin{equation} \label{eq:Del_ex}
  2 \pi i \sqrt{(\alpha + \beta \lambda^2)
  (\beta + \alpha \lambda^2)} =
  q - \bar{q}\lambda^2 + \lambda f_q,
\end{equation}
which in turn yields 
$a^2 = \sqrt{(\alpha + \beta \lambda^2)/
(\beta + \alpha \lambda^2)}$. 
Evaluating \eqref{eq:Del_ex} at $\lambda = 0$ gives 
$q = 2\pi i \sqrt{\alpha\,\beta}$. 
If $\rho:= \min \{ |\sqrt{\alpha/\beta}| ,\, 
| \sqrt{\beta/\alpha}|\}$, 
then $h \in \LGP$ for $0<r<\rho$. 
It is easily verified that 
$p(\lambda) = 2 \pi i \sqrt{- \det D}$ 
vanishes at the singularities of $h$. 
In summary, the matrix that dresses the vacuum 
into a Delaunay surface with 
neck radius $\omega$ is given by 
\begin{equation}
  h(\lambda) = \mathrm{diag}[\, 
  \sqrt[4]{\tfrac{\alpha + \beta \lambda^2}
  {\beta + \alpha \lambda^2}},\, 
  \sqrt[4]{\tfrac{\beta + \alpha \lambda^2}
  {\alpha + \beta \lambda^2}}\, ] 
\end{equation} 
Untwisting $h$, factoring 
$h = \mathcal{M}\mathcal{C}$ as 
in \eqref{eq:diagonal_h-split} and twisting back 
gives
\begin{equation} 
  \mathcal{M} = \frac{1}{2} \begin{pmatrix} 
  1+ \sqrt{\tfrac{\alpha + \beta \lambda^2}
  {\beta + \alpha \lambda^2}} & 
  \lambda - \lambda \sqrt{\tfrac{\alpha + 
  \beta \lambda^2}{\beta + \alpha \lambda^2}} \\ 
  \lambda^{-1}\sqrt{\tfrac{\beta + 
  \alpha \lambda^2}{\alpha + \beta \lambda^2}} 
  -\lambda^{-1} & 
  1 + \sqrt{\tfrac{\beta + \alpha \lambda^2}
  {\alpha + \beta \lambda^2}} \end{pmatrix}.
\end{equation}
Evidently $\mathcal{M}^* = \mathcal{M}^{-1}$ 
so that this is actually the decomposition 
according to Theorem \ref{th:h-split_cyl}. 
The domain of $\mathcal{M}$ is the genus one 
hyperelliptic curve $\mu^2  + \det D = 0$, 
the \textit{spectral curve} 
of the underlying Delaunay surface. 


\newsection{\bf Finite Blaschke Products.} 
The function $f_q$ occurring in Theorem 
\ref{th:f_q} is not always as easy to deal with as 
in the above example nor can it generally be 
explicitly computed. Nonetheless, even the 
trivial case $f_q \equiv 0$ is quite interesting 
and the corresponding dressing matrices are of such 
a simple form that the involved Iwasawa 
decomposition can be explicitly computed.  

\mylemma{}\label{th:f_q=0} 
Let $h \in \LGP$ and $F_c$ be the unitary frame 
of the round cylinder with monodromy $\chi_c(\tau_q) 
= \exp((q\lambda^{-1}-\bar{q}\lambda)A)$ for some 
$q \in \C^*$ with 
\begin{equation} \label{eq:q_condition}
  \sqrt{\bar{q}} - \sqrt{q} \notin \pi \,\Z.
\end{equation}
Assume $F = h \# F_c$ has a monodromy 
$h\,\chi_c(\tau_q)\,\exp(f_q\,A)\,h^{-1} \in \Lhol$ 
as in equation \eqref{eq:dress_monodromy}. 
Then the following are equivalent 
\begin{enumerate}
\item $h \chi_c(\tau_q) h^{-1} \in \Lhol$ 
\item $f_q \equiv 0$
\item $h\,A\,h^{-1}$ is rational on $\CP$ with only 
      simple poles.
\end{enumerate}

\Proof
We first show that (i) implies (ii). 
The assumption (i) allows us to write the 
monodromy of $F$ as 
$h\,\chi_c(\tau_q)\,h^{-1} \cdot 
h\exp(f_q\,A)\,h^{-1}$, where the first group 
of factors is unitary and the second group of 
factors is contained in $\LGP$ and has first 
term $I$. Therefore the condition that the whole 
expression is unitary is equivalent to 
$h\exp(f_q\,A)\,h^{-1}$ being unitary. 
But it is also in $\LGP$. Therefore it is in $\U1$.
But since it starts with $I$, it is identically 
equal to $I$. As a consequence $\exp(f_q\,A) = I$.
Hence $f_q \equiv 2 \pi k$ for some $k \in \Z$ 
and $k =0$, since $f_q$ is an odd function by 
Lemma \ref{th:f_q}. This shows that (i) implies (ii).

Let us prove that (ii) implies (iii). Recall from 
\eqref{eq:sinhcosh} and \eqref{eq:alphabeta} that 
$h\,\exp(p_q\,A)\,h^{-1} = \alpha\,\Id + 
\beta\,h\,A\,h^{-1}$. Since $f_q =0$ we have 
$p_q = \lambda^{-1}q - \lambda \bar{q}$ and 
therefore 
$\beta = \sinh(\lambda^{-1}q-\lambda\bar{q})$ is 
holomorphic on $\C^*$. 
By assumption, $\beta\,h\,A\,h^{-1}$ 
is holomorphic on $\C^*$ so that 
$h\,A\,h^{-1}$ is meromorphic on $\C^*$. 
Let $\lambda_0 \in \C^*$ be a root of 
$\beta = \sinh(\lambda^{-1}q-\lambda\bar{q})$, 
that is, it lies in the set
\begin{equation} \label{eq:spectral_set}
  \mathcal{S}_r^q = \left\{ \lambda \in A_r : 
  \lambda^{-1}q - \lambda\,\bar{q} \in \pi i \Z 
  \right\}. 
\end{equation}
Then $\beta'(\lambda_0) = \pm(\lambda_0^{-2}q - 
\bar{q}) \neq 0$ if and only if 
$\lambda_0 \neq \pm i \sqrt{q/\bar{q}}$. Our 
assumption \eqref{eq:q_condition} ensures that 
$\pm i \sqrt{q/\bar{q}} \notin \mathcal{S}_r^q$. 
Hence all roots 
of $\beta$ in $\C^*$ are simple. 
Therefore the entries of $h\,A\,h^{-1}$ 
can only have simple poles in $\C^*$ which 
must lie in $\mathcal{S}_r^q$.
Note also that $h\,A\,h^{-1}$ is holomorphic at 
$\lambda = 0$. Hence $h\,A\,h^{-1}$ is a holomorphic 
germ at the origin with a meromorphic 
extension to $\C^*$. In Section 3.5 of 
\cite{DorH:per} it is shown that the squares of the 
entries of $h\,A\,h^{-1}$ are finite at $\infty$. 
Hence $h\,A\,h^{-1}$ is rational on $\CP$ and 
proves that (ii) implies (iii).

Finally we show that (iii) implies (i). 
By equations \eqref{eq:sinhcosh} and 
\eqref{eq:alphabeta} we write 
$\chi = \alpha\,\Id + \beta\,hAh^{-1}$ 
and use the fact that 
$\beta\,h\,A\,h^{-1}$ is holomorphic on $\C^*$. 
With assumption (iii) this implies that $\beta$ 
is meromorphic on $\C^*$.
Since $\alpha$ is holomorphic 
on $\C^*$, so is $\alpha^2$ and consequently 
$\beta^2 = \alpha^2 - 1$ is holomorphic on 
$\C^*$, therefore $\beta$ itself is holomorphic 
on $\C^*$. Further, from 
$\chi^* = \chi^{-1} = \alpha\,\Id - 
\beta\,hAh^{-1}$ we deduce that $\alpha$ is real 
on $S^1$, i.e. $\alpha^* = \alpha$. 
And since $\beta^2$ is real and non-positive 
on $S^1$, and since $\beta$ is holomorphic on 
$\C^*$, we obtain $\beta^* = - \beta$. Consequently, 
$(hAh^{-1})^* = hAh^{-1}$ and thus 
$h\,\chi_c(\tau_q)\,h^{-1} = 
\exp((\lambda^{-1}q - \lambda\,q)h\,A\,h^{-1}) 
\in \Lhol$. This proves that (iii) implies (i) 
and concludes the proof of the lemma. \QED

\myremark{}
1. The result above seems to indicate that for 
every "type of $f_q$ " one obtains an associated 
type of dressing matrix.

2. In the case $f_q = 0$ and the setting of 
\cite{DorH:per}, the hyperelliptic surface 
associated with the dressing is simply the 
Riemann sphere $S^2$. 
If this class would contain some torus, then one 
would have an example for the 
"singular tori question" of \cite{Bob:tor}.


3. The set $\mathcal{S}_r^q$ in 
\eqref{eq:spectral_set} is finite for fixed 
$r$. Further $\mathcal{S}_r^q \subset \R$ if and 
only if $q \in i\R$, in which case $q$ also 
satisfies \eqref{eq:q_condition}. Finally, 
it is easy to see that if 
$\lambda \in \mathcal{S}_r^q$ then also 
$\bar{\lambda},\,1/\lambda,\,1/\bar{\lambda} 
\in \mathcal{S}_r^q$, and in particular, 
$\left(\mathcal{S}_r^q \right)^* = \mathcal{S}_r^q$.

We apply Lemma \ref{th:f_q=0} to a class of 
diagonal dressing matrices for which $f_q \equiv 0$.

\mycorollary{} 
Let $h = \mathrm{diag}[a,\,1/a] \in \LGP$. 
Then $h$ has an extension to $A_{r'}$ for some 
$r' \in [r,\,1]$ such that $h^* = h^{-1}$ on 
$A_{r'}$ and $h \chi_c(\tau_q) h^{-1} \in \Lhol$ 
if and only if $a^2$ is a finite Blascke product  
\begin{equation}\label{eq:a^2}
  a^2 (\lambda) = \prod_{j=1}^N 
  \frac{\alpha_j^2 - \lambda^2}
  {1 - \bar{\alpha}_j^2 \lambda^2},\,\alpha_j 
  \in \mathcal{S}_r^q. 
\end{equation}
\Proof Assume $h$ has a unitary branch on $S^1$ 
and $h \chi_c(\tau_q) h^{-1} \in \Lhol$. Then by 
Lemma \ref{th:f_q=0} we conclude that
$h\,A\,h^{-1} = \mathrm{off}[a^2,\,1/a^2]$ 
is rational with simple poles and zeroes located in 
$\mathcal{S}_r^q$. Hence $a^2$ is of the form 
\eqref{eq:a^2}. The converse is proven by direct 
verification.
\QED

Dressing matrices characterised in 
the previous proposition are thus of the form
\begin{equation} \label{eq:prod_simple}
  h = \prod_{j=1}^N \begin{pmatrix}  
  \sqrt{\tfrac{\alpha_j^2 - \lambda^2}
  {1 - \bar{\alpha}_j^2 \lambda^2}} & 0 \\ 0 & 
  \sqrt{\tfrac{1 - \bar{\alpha}_j^2 \lambda^2}
  {\alpha_j^2 - \lambda^2}} \end{pmatrix},\,
  \alpha_j \in \mathcal{S}_r^q
\end{equation}
and are a special instance of an interesting 
class to which we now turn our attention. 


\newsection{\bf Simple factors, untwisted case.} 

In this section we will deal primarily with 
untwisted loops. The previous discussion has 
naturally lead us to a class of special dressing 
matrices, the so called \textit{simple factors} of 
Uhlenbeck \cite{Uhl} and discussed in similar 
context in \cite{TerU} and \cite{Bur:iso}. 
The main feature is that dressing with 
simple factors is explicit. 
This construction, due to Terng and  
Uhlenbeck \cite{TerU} goes as follows:

\par We decompose 
$\C^2 = \mathrm{L} \oplus \mathrm{L}^\perp$ 
for $\mathrm{L} = \C \left( \begin{smallmatrix} 
a \\ b \end{smallmatrix} 
\right)$. Then for all $A \in \GL$ we have 
\begin{equation} \label{eq:perp}
  \bar{A}^t \mathrm{L} \perp 
  A^{-1} \mathrm{L}^\perp
\end{equation}
The hermitian projection 
$\pi_\l:\C^2 \to \mathrm{L}$ onto $\mathrm{L}$ is 
given by 
\begin{equation*}
 \pi_\l =  \frac{1}{|a|^2+|b|^2} \begin{pmatrix}
  |a|^2 & a\bar{b} \\ \bar{a}b & |b|^2 
  \end{pmatrix}.
\end{equation*}
Let $\alpha \in I_1^* = \left\{ \lambda \in \C : 
0 < |\lambda | < 1 \right\}$. Then the map 
\begin{equation*}
  \tau_\alpha(\lambda) = \frac{\alpha - \lambda}
  {1-\bar{\alpha}\lambda}
\end{equation*}
is invariant under $\varrho$. For given 
$\mathrm{L}$ and $\alpha \in I_1$, a 
\textit{simple factor} \cite{TerU} is a 
loop of the form 
\begin{equation} \label{eq:psi}
  \psi_{\alpha,\l}(\lambda) = \pi_\l + 
  \tau_\alpha(\lambda) \,\pi_\l^\perp.
\end{equation}
By construction, $\psi_{\alpha,\l}: \CP \setminus 
\left\{ \alpha,\,1/\bar{\alpha} \right\} \to \GL$ 
is analytic and since $|\alpha |>0$, 
\begin{equation}
  \psi_{\alpha,\l} \in \Lambda_r^+ \GL
  \mbox{ for } r < | \alpha |.
\end{equation}
Clearly, simple factors are not twisted.
Further, since $\psi_{\alpha,\l}^* = 
  \psi_{\alpha,\l}^{-1} = 
  \pi_\l + \tau_\alpha^{-1} \pi_\l^\perp$, 
we have that 
\begin{equation}
  \psi_{\alpha,\l} \in \Lambda_r \mathbf{U}(2)
  \mbox{ for } r > | \alpha |.
\end{equation}
For later use we also 
note the fact that for any $A \in \Utwo$ 
we have 
\begin{equation} \label{eq:change}
  \psi_{\alpha,A\l} = A \,\psi_{\alpha,\l}\,A^{-1}.
\end{equation}
Let $U \in  \Lambda_{\mbox{\tiny{$\C^*$}}}\SU,\,
\alpha \in I_1^*$ 
and $\mathrm{L} \in \C\mathbb{P}^2$ and  
$\psi_{\alpha,\l}$ be the corresponding simple 
factor. Then it can be shown, see e.g \cite{KilSS}, 
that 
\begin{equation}
  \psi_{\alpha,\l}\,U\,
  \psi_{\alpha,\overline{U(\alpha)}^t \l}^{-1} \in 
  \Lambda_{\mbox{\tiny{$\C^*$}}}\SU
\end{equation}
Consequently, we have an explicit 
$r$--Iwasawa decomposition of $\psi_{\alpha,\l}\,U$ 
for $r < |\alpha |$, given by 
\begin{equation} \label{sdf}
  \psi_{\alpha,\l}\,U = \left( \psi_{\alpha,\l}\,U\,
  \psi_{\alpha,\overline{U(\alpha)}^t \l}^{-1} 
  \right) \psi_{\alpha,\overline{U(\alpha)}^t \l}.
\end{equation}
For the geometric applications considered in 
this paper it is of great value to have
a simple and  explicit formula for the dressed 
frame. Simple factors, for which $\mathrm{L}$ 
depends on $\lambda$ and are of the form 
\begin{equation}
  \psi_{\alpha,\l(\lambda)} =  \pi_{\l(\lambda)} + 
  \tau_\alpha(\lambda) \,(\Id - \pi_{\l(\lambda)})
\end{equation}
will be called \emph{generalised simple factors}. 
We would like to present an analogous result 
to \ref{sdf} for such generalised simple factors: 

\mytheorem{}
Let $\mathrm{L}: \C^* \rightarrow \CP$ be holomorphic and 
$\langle \mathrm{L},\,\mathrm{L}^* \rangle \neq 0$ for all 
$\lambda \in \C^*$. 
Let $\alpha \in I_1^*$ and $\psi_{\alpha,\l(\lambda)}$ 
be the corresponding generalised simple factor. 
Then for any $U \in  \Lambda_{\mbox{\tiny{$\C^*$}}}\SU$
\begin{equation} \label{gsdf}
\psi_{\alpha,\l(\lambda)}\,\# \,U = 
  \psi_{\alpha,\l(\lambda)}\,U\,
  \psi_{\alpha,\overline{U(\alpha)}^t 
  \l(\alpha)}^{-1}. 
\end{equation}
\Proof
Fix $\lambda_0 \in \C^*$ and write 
$\mathrm{L}(\lambda) = [a(\lambda) : b(\lambda)]$. Then for 
$\mathrm{L}_0 = [1 : 0]$ and the holomorphic map
\begin{equation}
  W(\lambda) = \begin{pmatrix} 
  a(\lambda) & -\overline{b(1/\bar{\lambda})} \\
  b(\lambda) & \overline{a(1/\bar{\lambda})} \end{pmatrix}
\end{equation}
we have $\mathrm{L}(\lambda) = W(\lambda)\,\mathrm{L}_0$ and 
$W(\lambda) \in \Lambda_{\mbox{\tiny{$\C^*$}}} \Utwo$. 
Applying \ref{eq:change} to 
$A = W(\lambda), \,\lambda \in \C^*$, we obtain 
$\psi_{\alpha,\l(\lambda)} = 
W(\lambda)\,\psi_{\alpha,\l_0}\,W(\lambda)^{-1}$ 
and 
\begin{align*}
  \psi_{\alpha,\l(\lambda)}\,U(\lambda) &= 
  W(\lambda)\psi_{\alpha,\l_0}W(\lambda)^{-1} 
  U(\lambda) \\ &= \left( W(\lambda)\, 
  \psi_{\alpha,\l_0}\, W(\lambda)^{-1} U(\lambda)\, 
  \psi_{\alpha,\l_0'}^{-1} \right) 
  \psi_{\alpha,\l_0'},
\end{align*}
where $\mathrm{L}_0' = \overline{(W(\alpha)^{-1} 
U(\alpha) )}^t\mathrm{L}_0$.
The right side is an Iwasawa decomposition. 
The unitary part can be rewritten in the form
\begin{equation*}
  W(\lambda)\,\psi_{\alpha,L_0}\,W(\lambda)^{-1} 
  \,U(\lambda)\,\psi_{\alpha,L_0'}^{-1} =
  \psi_{\alpha,L(\lambda)}\,U(\lambda)\,
  \psi_{\alpha,L_0'}^{-1}. 
\end{equation*}
It is straightforward to verify 
$\mathrm{L}_0' = \overline{U(\alpha)}^t 
\mathrm{L}(\alpha)$, thus proving \eqref{gsdf}. \QED
%


\newsection{\bf Twisting of simple factors}

Since generalised simple factors are in general 
untwisted, we need to
modify the concept so that it can also be used for 
the twisted case. We discuss two approaches to 
deal with simple factors in the twisted case by 
firstly, taking a simple factor and applying the 
'twisting map' and secondly, looking for products 
of simple factors which happen to be twisted. 
Applying the 'twisting map' \eqref{eq:twist} 
to a simple factor we obtain 
\begin{equation*}
  \psi_{\alpha,\l(\lambda)} =D(\lambda)\,
  \psi_{\alpha,\l (\lambda^2)}\,D(\lambda)^{-1} =  
  \psi_{\alpha, D(\lambda) \l (\lambda^2)}.
\end{equation*}
\mylemma{}
a) Twisting a generalised simple factor 
produces again a generalised simple factor.

b) Let $\psi$ be a simple factor (L is 
independent of $\lambda$). Then the 
corresponding twisted simple
factor is again independent of 
$\lambda$ if and only if
\begin{equation}\label{eq:picase} 
  D\,\pi_\l \,D^{-1} = \pi_\l \mbox{ or } 
  \pi_\l^\perp. 
\end{equation}
\QED

\myremark{}
Since $D=D(\lambda)$ is diagonal, the only 
projections satisfying the condition 
\eqref{eq:picase} above are 
those onto the canonical basis vectors.


\newsection{\bf Twisted products of simple factors}
As mentioned earlier we also consider products of 
two simple factors and determine, when such 
matrices are twisted. Recall the two 
involutions $\varrho,\,\s$ defined in 
\eqref{eq:sigma} respectively \eqref{eq:varrho}. 

\myprop{} \label{th:twist}

a) Let $g \in \Lambda_r \SL$. 
Then $(\s g) g$ is twisted if and 
only if $[\,\s g,\,g\,]=0$.

b) If $g \in \LGU$, then $(\s g)g \in \LGU$.

c) Let $\psi_{\alpha,\l}$ and 
$\psi_{\beta,\widehat{\l}}$ be simple factors.
If $\sigma (\psi_{\alpha,\l})\, 
\psi_{\beta,\widehat{\l}}$ is twisted, then 
$\psi_{\alpha,\l} = \psi_{\beta,\widehat{\l}}$.

d) If $\sigma (\psi_{\alpha,\l})\psi_{\alpha,\l}$ 
is twisted, then either $\sigma(\pi_\l) = \pi_\l$ 
or $\sigma(\pi_\l) = \pi_\l^\perp$.

\Proof a) The loop $(\s g) g$ is twisted if and 
only if $\s \left( (\s g)\,g \right) = g \,\s g = 
(\s g)\,g$, which is equivalent to 
$[\s g,\,g] = 0$. 

b) The loop $g \in \LGU$ if and only if 
$\varrho \,g = g$. Then 
$\varrho \left( (\s g)g \right) = 
(\varrho \s g)\,(\varrho g) = 
(\s \,\varrho \,g ) \, (\varrho\, g) = (\s g)\,g$, 
since $[\varrho,\,\sigma ] = 0$. Hence 
$(\s g)\,g \in \LGU$. 

c) The product is twisted if and only if 
$\sigma (\psi_{\alpha,L}) \cdot  
\psi_{\beta,\widehat{L}} =  \psi_{\alpha,L} 
\cdot \sigma ( \psi_{\beta,\widehat{L}})$. 
Let us abbreviate $A = \psi_{\alpha,\l}$ and  
$B = \psi_{\beta,\widehat{\l}}$. 
Then the pole of $A$ is at 
$1/\overline{\alpha}$, while the pole of 
$B$ is at  $1/\overline{\beta}$. 
On the other hand, the pole of $\sigma (A)$ 
is at $- 1/\overline{\alpha}$ and the 
analogous result holds for $B$.
Thus comparing the two sides we obtain 
$\alpha = \beta$. Expanding near a pole 
and comparing the factors we derive $\mathrm{L} = 
\widehat{\mathrm{L}}$.

(d) If $(\s \psi_{\alpha,\l} )\psi_{\alpha,\l}$ 
is twisted, then by part (a), 
$[\s \psi_{\alpha,\l},\,\psi_{\alpha,\l}]=0$ 
which is equivalent to the eigenspaces of 
$\s \psi_{\alpha,\l}$ and 
$\psi_{\alpha,\l}$ coinciding. Thus there are two 
possibilities:
\begin{equation}\label{eq:picase1} 
  \s_3\,\pi_\l \,\s_3^{-1} = \pi_\l \mbox{ or } 
  \pi_\l^\perp. 
\end{equation}
This proves (d) and concludes the proof of the 
proposition. \QED

We are going to evaluate the two possibilities 
in \eqref{eq:picase1}. 

(i) Let us turn to the first case in 
\eqref{eq:picase1}: 
For the line 
$\mathrm{L} = \C \bigl( \begin{smallmatrix} 
1\\0 \end{smallmatrix} \bigr)$ we have 
$\psi_{\alpha,\l}= \mathrm{diag}[1,\,\tau_\alpha]$. 
Using $\tau_\alpha(\lambda)\tau_\alpha(-\lambda) =  
  \tau_{\alpha^2}(\lambda^2)$ 
we obtain $(\s \psi_{\alpha,\l})\psi_{\alpha,\l} = 
\pi_\l + \tau_{\alpha^2}(\lambda^2)\pi^\perp_\l$. 
Dividing by the square root of the determinant  
we arrive at all diagonal twisted simple factors 
\begin{equation} \label{eq:simple1}
 g_{\alpha,\l}(\lambda) = 
  \sqrt{\tau_{\alpha^2}^{-1}
  (\lambda^2)} \, \pi_\l + \sqrt{\tau_{\alpha^2}
  (\lambda^2)}\, \pi^\perp_\l.
\end{equation}
Then $g_{\alpha,\l}(\lambda) \in \LGP$ for 
$0< r< |\alpha|$ and in matrix form is given by 
\begin{equation}\label{eq:simple1_matrix}
  g_{\alpha,\l}(\lambda) = \begin{pmatrix} 
  \sqrt{\tfrac{1-\bar{\alpha}^2\lambda^2}
  {\alpha^2 - \lambda^2}} & 0 \\ 
  0 & \sqrt{\tfrac{\alpha^2 - \lambda^2}
  {1-\bar{\alpha}^2\lambda^2}}
  \end{pmatrix} 
\end{equation}
Notice that $g_{\alpha,\l}(\lambda)^{-1}$ 
corresponds to a factor of the matrix 
given in \eqref{eq:prod_simple}.

(ii) Let us turn to the second case in 
\eqref{eq:picase1}: $\s_3 \pi_\l \s_3^{-1} = 
\pi_\l^\perp$. 
For $\mathrm{L} =\C ( a,\,b )^t$ and hermitian 
projection $\pi_\l:\C^2 \to \mathrm{L}$ we have 
$\sigma \pi_\l = \pi_\l^\perp$ if and 
only if $|a|^2=|b|^2 = 1/2$. 
Setting $a=\tfrac{1}{\sqrt{2}}\exp(is)$ and 
$b=\tfrac{1}{\sqrt{2}}\exp(it)$ and 
rescaling, we may assume without loss of 
generality that $\mathrm{L} = 
\C \bigl( e^{i\theta}, \, 1 \bigr)^t$. Then 
$\pi_\l = \tfrac{1}{2}(\Id + \mathrm{off}
[ e^{i\theta},\,
e^{-i\theta}])$ and for 
$\psi_{\alpha,\l} = \pi_\l + \tau_\alpha(\lambda)\, 
\pi_\l^\perp$ 
we have $\sigma \psi_{\alpha,\l} = \pi_\l^\perp 
+ \tau_\alpha(-\lambda)\,\pi_\l$ and consequently 
\begin{align*}
  (\sigma \psi_{\alpha,\l})\psi_{\alpha,\l} &= 
  \tau_\alpha(\lambda)\,\pi_\l^\perp + 
  \tau_\alpha(-\lambda)\,\pi_\l \\
  &= \tfrac{1}{1-\bar{\alpha}^2\lambda^2} (
  (\alpha - \bar{\alpha}\lambda^2)\,\Id
  + \lambda (1 - |\alpha|^2)\,
  \mathrm{off}[e^{i\theta},\,e^{-i\theta}]).
\end{align*}
Normalising, we arrive at the second class of 
twisted simple factors in $\LGP$ for 
$0 < r < |\alpha|$, given by
\begin{equation} \label{eq:simple2}
 g_{\alpha,\l}(\lambda) = 
  \sqrt{\tau_\alpha^{-1}(\lambda)
  \,\tau_\alpha (-\lambda)} \, \pi_\l + 
  \sqrt{\tau_\alpha (\lambda) 
  \,\tau_\alpha^{-1} (-\lambda)}\, \pi^\perp_\l.
\end{equation}
In matrix form, these simple factors look like
\begin{equation} \label{eq:simple2_matrix}
  g_{\alpha,\l}(\lambda) = 
  \sqrt{\frac{1-\bar{\alpha}^2\lambda^2}
  {\alpha^2 - \lambda^2}}
  \begin{pmatrix} \alpha - \bar{\alpha}\lambda^2 & 
  \lambda (1- |\alpha|^2)e^{i\theta} \\ 
  \lambda (1 - |\alpha|^2)e^{-i\theta} & 
  \alpha - \bar{\alpha}\lambda^2 
  \end{pmatrix}.
\end{equation}
The simple factors $g_{\alpha,\l}$ derived 
in \eqref{eq:simple1} and \eqref{eq:simple2} 
are uniquely determined by their 'singularity' 
$\alpha \in I_1$ and choice of line 
$\mathrm{L} \in \CP$. Note also that in both cases 
\eqref{eq:simple1} and \eqref{eq:simple2} we have 
that 
\begin{equation} \label{eq:simpleunitary}
  g_{\alpha,\l} \in \LGU 
  \mbox{ for } |\alpha| < r \leq 1.
\end{equation}
The key aspect is that dressing with simple factors 
is explicit. This idea is due to 
Terng and Uhlenbeck \cite{TerU} and has lead to 
variants as in \cite{Bur:iso} and \cite{KilSS}. 
The version we need is proven by 
twisting Theorem 1.2 in \cite{KilSS} and is as 
follows 

\mytheorem{} Let $M$ be a connected Riemann 
surface with universal cover $\MT$ and let 
$F(z,\lambda) \in \mathcal{F}(\MT)$. 
Let $g_{\alpha,\l} \in \LGP$ be a 
simple factor. Then 
\begin{equation} \label{eq:simple_dress}
  g_{\alpha,\l} \# F = 
  g_{\alpha,\l} \,F \,g^{-1}_{\alpha,\l^\prime}
\end{equation}
where $\mathrm{L}^\prime = 
\overline{F(z,\alpha)}^t \mathrm{L}$ and 
$g^{-1}_{\alpha,\l^\prime}$ is, pointwise in 
$z \in \MT$, a simple factor of the same 
form as $g_{\alpha,\l}$. \QED

We use this factorisation theorem to characterise 
when it is possible to dress an extended unitary 
frame with trivial isotropy by a simple factor 
while retaining a symmetry.

\mycorollary{} \label{th:invariance} 
Let $\chi(\tau) \in \Lhol$ be a monodromy 
of an extended unitary frame 
$F \in \mathcal{F}(\MT)$. 
Let $g_{\alpha,\l} \in \LGP$ be a simple factor. 
Then 
$g_{\alpha,\l}\,\chi(\tau) \,g_{\alpha,\l}^{-1}$ 
is a monodromy of $g_{\alpha,\l} \# F$ 
with respect to $\tau$ if and only if 
\begin{equation} \label{eq:invariance}
  \overline{\chi(\alpha,\tau)}^t 
  \mathrm{L} = \mathrm{L}.
\end{equation}
\Proof Let $\hat F = g_{\alpha,\l} \# F$. 
Then $\hat F = g_{\alpha,\l}\,F\,
g^{-1}_{\alpha,\l^\prime}$ with 
$\mathrm{L}^\prime = \overline{F(z,\alpha)}^t 
\mathrm{L}$ by \eqref{eq:simple_dress}. Let 
$\tau^* F = \chi(\tau)\,F\,k$. Then 
\begin{align*}
  \tau^* \hat F &= g_{\alpha,\l}\,(\tau^* F)\, 
  g^{-1}_{\alpha,\tau^*\l^\prime} \\ 
	   &= g_{\alpha,\l}\,\chi(\tau) \,F\, 
  k\,g^{-1}_{\alpha,\l''} \mbox{ with } 
  \mathrm{L}'' = \bar{k}^t\,
  \overline{F(\alpha)}^t \,
  \overline{\chi(\alpha,\tau)}^t \mathrm{L}  = 
  \tau^* \mathrm{L}' \\
  &= g_{\alpha,\l}\,\chi(\tau)\, g_{\alpha,\l}^{-1}
  \,g_{\alpha,\l} \,F \, g^{-1}_{\alpha,\l_1 } 
  k \mbox{ by } \eqref{eq:change}
\end{align*}
where $\mathrm{L}_1 = \overline{F(z,\alpha)}^t\,
\overline{\chi(\alpha,\tau)}^t \mathrm{L}$. 

If $g_{\alpha,\l}\,\chi(\tau)\, g_{\alpha,\l}^{-1}$ 
is a monodromy of $\hat F$ and we write 
$\tau^* \hat{F} = g_{\alpha,\l}\,\chi(\tau)\, 
g_{\alpha,\l}^{-1} \, \hat{F} \,\hat{k}$ 
and combine with the above, we obtain 
$g_{\alpha,\l_1}^{-1}\,k = 
  g_{\alpha,\l'}^{-1} \,\hat{k}$. 
Hence $g_{\alpha,\l_1}g_{\alpha,\l'}^{-1}$ 
is independent of $\lambda$ and since there exists a 
$\lambda_0 \in \CP$ at which $\tau_\alpha(\lambda_0) = 1$ we obtain 
$g_{\alpha,\l_1}g_{\alpha,\l'}^{-1} = \Id$.  
This implies that $L' = L_1$ which implies \eqref{eq:invariance}. 
The converse is proven by direct verification. 
\QED


\section{Dressing non--finite type surfaces}

In this section we discuss a family of CMC cylinders 
found in \cite{KilMS} that arise as perturbations 
of Delaunay surfaces. 
The resulting surfaces may possess an arbitrary 
number of umbilics and are thus not of finite type and  
are CMC cylinders with one Delaunay end \cite{KilKRS}, see Figure 1. 
We can dress this class of cylinders with simple factors, having the 
effect of adding bubbles to the surface, see Figure 2.

\begin{figure}\label{fig:pert}
  \centering
    \includegraphics[scale=1.5]{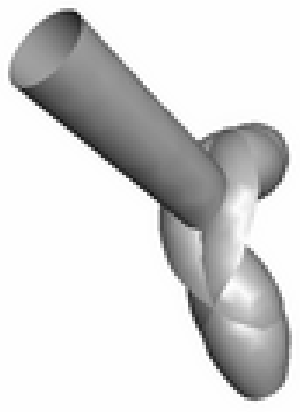}
    \includegraphics[scale=1.3]{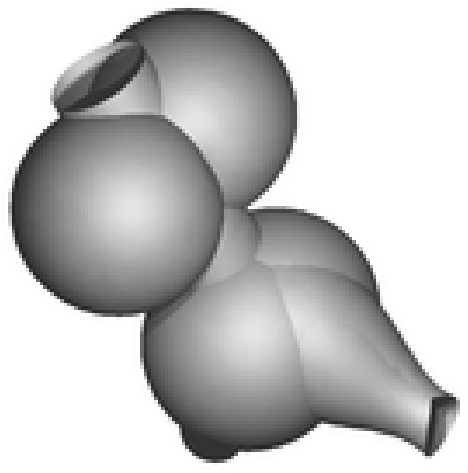}\\ 
    \caption{\footnotesize A perturbed round cylinder on the left. A perturbed 
   and dressed round cylinder on the right. 
   Images generated with \emph{CMCLab} \cite{Sch:cmclab}. For more images see 
  \cite{Sch:gallery}.}
\end{figure}

\newsection{\bf Dressing cylinders with umbilics.}
We modify a standard result from the theory of 
differential equations with regular singular points. 
Since the eigenvalues are $\lambda$--dependent, 
we can avoid the assumption that two elements in the 
spectrum not differ by an integer by working on an 
appropriate $\lambda$--circle $C_r$.

\mytheorem{} \label{th:zap}
Let $U_0^* \subset \C$ be an open neighbourhood of 
$z=0$ and $\xi \in \Lambda \Omega(U_0^*)$ with a 
simple pole at $z=0$ and residue
$\mathrm{res}_{\,0} \xi = D$ 
where $D$ has the form \eqref{eq:delaunay}.
Moreover, we assume that $D$ 
satisfies the first closing condition 
\eqref{eq:closing1}. Then there exists 
an $r \in (0,1)$ and a solution of 
$d\Psi = \Psi \xi$ of the form $\Psi = z^D P$ 
with $P:U_0 \to \LGC$ holomorphic. 
Further, $\Psi$ has Delaunay monodromy and is the 
holomorphic frame associated with $(\xi,\,P(1),\,1)$.

\Proof
Let $U_0 = U_0^* \cup \{ 0 \}$ and write 
$\xi = D \tfrac{dz}{z} + \eta$ with 
$\eta \in \Lambda \Omega_{U_0}$. The differential 
equation for $P$ is 
$dP = P \eta + [ P,\, D \tfrac{dz}{z}]$. 
We will show that this differential equation has a 
holomorphic solution. Expanding 
$P = \sum_{k = 0}^{\infty} P_k z^k$ and 
$\xi = (\tfrac{1}{z}D +\Sigma \eta_k z^k )dz$ 
at $z=0$, the coefficients are recursively given 
by 
$$
kP_k + [D , P_k ] = \sum_{r+s=k-1} P_r \eta_s.
$$
By our Ansatz $P_k = 0$ for $k < 0$ and we are free 
to choose $P_0 \in \LGC$ with $[D , P_0 ]=0$. 
(Note this freedom only means 
$P_0 = \alpha \Id + \beta D$, 
since $D$ is a semisimple $2 \times 2-$matrix
with different eigenvalues.)
The eigenvalues of $D$, are of the form $\pm \mu$, 
where $\mu$ is a solution to the equation 
$\mu^2 = - \det D$. Then 
the operators $k \rm{Id} + \rm{ad}_{D}$, 
$k \in \Z$, have spectrum 
$\sigma = \{ k,\, k \pm  2\mu \}$. 
We are only interested in $k \geq 1$, in which case 
$0 \notin \sigma$ if and only if 
$\mu \notin \tfrac{1}{2} \N$. 

We need to show that there exists an 
$r \in (0,\,1)$ for which 
$\mu (\lambda) \notin \frac{1}{2}\N$ for every 
$\lambda \in C_r$. 
A priori, $\mu(\lambda) \in \R$ if and only if 
either $\lambda \in S^1$ or $\lambda \in \R$. 
Since we seek $r \in (0,\,1)$ we have 
$\mu(\lambda) \in \R$ if and only if 
$\lambda \in \R$ we need to ensure 
$\mu(\pm r)\notin \tfrac{1}{2}\N$. 
Since $D$ satisfies the first closing 
condition by assumption, we know that 
$\mu(1) = n/2$ for some $n \in \Z$ 
the winding number of the Delaunay surface. 
A direct computation shows 
$\mu(\pm r) \notin \tfrac{1}{2}\Z$ 
if and only if 
$r^{-2} + r^2 \neq \pm (k^2 - n^2)/(4ab) + 2$ 
for all $k \in \N$. 
The sequences $a^{\pm}_k := 
\pm (k^2 - n^2)/(4ab) + 2$ 
are monotonic, hence we can always find an 
$0 < r <1$ for which $r^{-2} + r^2 
\notin \{ a^{\pm}_k \}$. 

This proves the existence of a formal power 
series solving the differential equation for $P$. 
The coefficients of this formal power series are 
defined for all $z \in U_0$. The domain of 
analyticity of $P$ is $U_0$ by standard 
ODE arguments, see e.g \cite{Har}. It is shown in 
\cite{DorW:noids} that $P$ has the twisted $\lambda$ 
behaviour and that the coefficients of $P$ as 
functions of $\lambda$ are in $\mathcal{A}_r$. 
To show that $\Psi$ has Delaunay 
monodromy, recall from \ref{sec:Delaunayexample} 
that the Delaunay monodromy is $\exp(2 \pi i D)$. 
For the translation $\tau(z) = z + 2\pi i$ 
the monodromy of $\Psi$ is $\tau^*\Psi \Psi^{-1} = 
\tau^*z^D \tau^*P P^{-1} z^{-D} = \exp(2 \pi i D)$.
\QED

For our purposes, the specific form of $\Psi$ 
proven above is of great importance to us, 
since it gives easy control over the monodromy 
singled out by the singularity at $D$.
In particular, if we dress the corresponding 
frame  by  some simple factor it is easy to
determine the dressed monodromy, at least if 
the associated surface has umbilical points. 
From a geometric point of view this is 
particularly interesting, since by a result 
of \cite{Mah}, the classical Bianchi--B\"acklund 
transformation does correspond to dressing 
by some simple factor. Actually, in the twisted 
setup, our 'simple factors' are in fact a product 
of two simple factors. These correspond to the 
two step procedure of the Bianchi--B\"acklund 
transformation. As a matter of fact,  
by making the  right choice for the singularity 
of the simple factor, is is possible to control the 
topology of the resulting surface. 

\myremark{}
1) The proof of the Theorem above actually shows 
that there is some $0 < r < 1$ such that for all 
$0 < r \leq r' <1$ the claim holds. 

2) Instead of choosing some $0 < r < 1$, 
one could ask for what potentials $\xi$ 
one can obtain a solution of the
form $ \Psi = z^D \, P$ without 
singularities on $S^1$. Once this is 
achieved, one can conclude as above and 
infer statements about the monodromy and 
thus one can clear the way for some 
understanding for the effect of dressing 
by simple factors. Such a condition has 
been given in \cite{DorW:noids}.

In order to dress a CMC surface with non-trivial 
topology into a new CMC surface with the same 
topology, we recall from \cite{DorH:per} when for a 
given triple $(\xi,\Phi_0,\tilde{z}_0)$ one 
specific member of the associated family is 
invariant under $\Delta$, thus characterising 
the period problem in the DPW framework: 

Let $f_\lambda$ be an associated family. 
Let $F$ be the extended unitary 
frame of the surface with monodromy 
$\chi : \Delta \to \Lhol$.
Then there exists a $\lambda_0 \in S^1$ such 
that $\tau^* f_{\lambda_0}= f_{\lambda_0}$ for all 
$\tau \in \Delta$ if and only if for all 
$\tau \in \Delta$, $\chi$ satisfies both
\begin{align}
\chi (\tau,\,\lambda_0) = 
\pm \rm{Id},&\label{eq:closing1}\\
\left. \tfrac{\partial}{\partial \lambda}
\right|_{\lambda_0} 
\chi(\tau,\,\lambda)
 = 0.& \label{eq:closing2}
\end{align}

\mycorollary{} 
Let $\xi \in \Lambda \Omega(\C^*)$ as in Theorem 
\ref{th:zap} and $g_{\alpha,\l} \in \LGP$ be a 
simple factor.  Choosing $r \in (0,\,1)$ as above 
and $\alpha \in \C$ with $|\alpha| \in (r,\,1)$ 
such that 
\begin{equation} \label{eq:det}
  \det D(\alpha) = \tfrac{n}{4} \mbox{ for some } n \in \Z,
\end{equation}
then $(\xi,\, g_{\alpha,\l}\, P(1), \,1)$ 
generates a CMC cylinder (with winding number $n$).

\Proof 
Let $\Phi = z^D\,P$ and $F$ be the holomorphic 
respectively unitary frame generated by 
$(\xi ,\, P(1),\, 1)$. 
Both have Delaunay monodromy 
$\chi (\tau,\,\lambda) = \exp(2\pi iD(\lambda))$ 
with respect to the translation 
$\tau(z) = z + 2\pi i$. 
The condition \eqref{eq:det} ensures that 
$\chi(\tau,\,\alpha) = \pm \Id$ and 
consequently \eqref{eq:invariance} holds 
for any choice of line $\mathrm{L} \in \CP$. 
Hence, by Corollary \ref{th:invariance}, 
$\hat\chi(\tau) = 
g_{\alpha,\l}\,\chi(\tau)\,g_{\alpha,\l}^{-1}$ 
is a monodromy of $g_{\alpha,\l} \# F$. 
The conditions \eqref{eq:closing1} and 
\eqref{eq:closing2} are verified using the facts 
$\chi(\tau,1) = \pm \mathrm{Id}$ and $\left. 
\partial_{\lambda}\chi(\tau)\right|_{\lambda=1}= 0$.
\QED


\newsection{\bf Dressing $3$-Noids.}
Let us briefly outline how the above theory can be 
applied to construct the dressed $3$-Noids of 
\cite{KilSS}. Let 
$\mathcal{T} = \CP \setminus \{0,\,1,\,\infty\}$ 
and $f: \mathcal{T} \to \R^3$ a CMC $3$-Noid with 
Delaunay ends as constructed in \cite{DorW:noids} 
or \cite{Sch:tri}. Let 
$F \in \mathcal{F}(\widetilde{\mathcal{T}})$ 
be an extended frame such that at $\lambda = 1$ the 
immersion $f$ is obtained via the Sym--Bobenko 
formula \eqref{eq:Sym}. 
Then there are three monodromies, 
$\chi_1,\, \chi_2,\,\chi_3$, one for each end, 
which can be computed in terms of 
$\Gamma$--functions \cite{DorW:noids}, 
\cite{KilSS} and satisfy 
$\prod_{j=1}^3 \chi_j = \Id$. 

In \cite{KilSS} it was shown that there exist 
values $\alpha \in \C^*,\,|\alpha |<1$ and 
invariant subspaces $\mathrm{L} \in \CP$ such that 
\eqref{eq:invariance} holds for all three 
$\chi_j$'s. Consequently, by Corollary 
\ref{th:invariance} the unitary frame 
obtained by dressing 
with the simple factor $g_{\alpha,\l}$ 
has monodromies $\hat \chi_j = 
g_{\alpha,\l} \, \chi_j \, g_{\alpha,\l}^{-1}$ and 
since the closing conditions \eqref{eq:closing1} and \eqref{eq:closing2} are invariant under conjugation, 
the resulting surface is again a 
CMC $3$--Noid.
%


\providecommand{\bysame}{\leavevmode\hbox to3em{\hrulefill}\thinspace}
\providecommand{\MR}{\relax\ifhmode\unskip\space\fi MR }
\providecommand{\MRhref}[2]{%
  \href{http://www.ams.org/mathscinet-getitem?mr=#1}{#2}
}
\providecommand{\href}[2]{#2}



\begin{thebibliography}{10}

\bibitem{Bob:tor}
A.~I. Bobenko, \emph{All constant mean curvature tori in {$\mathbb{R}^3$},
  {$\mathbb{S}^3$}, {$\mathbb{H}^3$} in terms of theta-functions}, Math. Ann.
  \textbf{290} (1991), 209--245.
\bibitem{BobPS}
A.~I. Bobenko, T.~V. Pavlyukevich and 
B. A. Springborn,
\emph{Hyperbolic constant mean curvature one 
surfaces: spinor representation and trinoids in 
hypergeometric functions}, Math.Z \textbf{245}, (2003), 63--91.


\bibitem{Bun}
L.~Bungart, \emph{On analytic fiber bundles}, Topology \textbf{7} (1968),
  55--68.

\bibitem{Bur:iso}
F.~Burstall, \emph{Isothermic surfaces: conformal geometry, {C}lifford algebras
  and integrable systems.}, Integrable systems, Geometry and Topology, vol.~80,
  taiwan, 1995, math.DG/0003096, pp.~353--382.

\bibitem{BurP:dre}
F.~Burstall and F.~Pedit, \emph{Dressing orbits of harmonic maps}, Duke Math.
  J. \textbf{80} (1995), 353--382.

\bibitem{DorH:per}
J.~Dorfmeister and G.~Haak, \emph{On constant mean curvature surfaces with
  periodic metric}, Pacific J. Math. \textbf{182} (1998), 229--287.

\bibitem{DorH:dre}
\bysame, \emph{Investigation and application of the dressing action on surfaces
  of constant mean curvature}, Q. J. Math. \textbf{51} (2000), 57--73.

\bibitem{DorH:sym2}
\bysame, \emph{On symmetries of constant mean curvature surfaces. {II}.
  symmetries in a {W}eierstrass-type representation}, Int. J. Math. Game Theory
  and Algebra \textbf{10} (2000), 121--146.

\bibitem{DorH:cyl}
\bysame, \emph{Construction of non-simply connected {CMC} surfaces via
  dressing}, J. Math. Soc. Japan \textbf{55} (2003), no.~2, 335--364.

\bibitem{DorPW}
J.~Dorfmeister, F.~Pedit, and H.~Wu, \emph{Weierstrass type representation of
  harmonic maps into symmetric spaces}, Comm. Anal. Geom. \textbf{6} (1998),
  no.~4, 633--668.

\bibitem{DorW:noids}
J.~Dorfmeister and H.~Wu, \emph{Construction of constant mean curvature
  $n$-noids from holomorphic potentials}, in preparation, 2002.

\bibitem{Har}
P.~Hartman, \emph{Ordinary differential equations}, John Wiley and Sons, Inc.,
  New York, 1964.

\bibitem{Kil:thesis}
M.~Kilian, \emph{Constant mean curvature cylinders}, Ph.D. thesis, Univ. of
  Massachusetts, Amherst, 2000.

\bibitem{Kil:del}
\bysame, \emph{On the associated family of {D}elaunay surfaces},
  Proc. AMS, in press.

\bibitem{KilKRS}
M.~Kilian, S.~Kobayashi, W.~Rossman, and N.~Schmitt, \emph{{CMC} surfaces in
  three dimensional space forms}, math.DG/0403366.

\bibitem{KilMS}
M.~Kilian, I.~McIntosh, and N.~Schmitt, \emph{New constant mean curvature
  surfaces}, Experiment. Math. \textbf{9} (2000), no.~4, 595--611.

\bibitem{KilSS}
M.~Kilian, N.~Schmitt, and I.~Sterling, \emph{Dressing {CMC} n-{N}oids}, 
  Math. Z., \textbf{246}, (2004), 501--519.

\bibitem{Kob}
S.~Kobayashi, \emph{Bubbletons in 3-dimensional space forms via the {DPW}
  method}, Master's thesis, Kobe University, 2001.

\bibitem{Mah}
A.~Mahler, \emph{{B}ianchi-{B}\"acklund and dressing transformations on
  constant mean curvature surfaces}, Ph.D. thesis, Univ. of Toledo, 2002.

\bibitem{McI}
I.~McIntosh, \emph{Global solutions of the elliptic 2d periodic {T}oda
  lattice}, Nonlinearity \textbf{7} (1994), no.~1, 85--108.

\bibitem{PinS}
U.~Pinkall and I.~Sterling, 
\emph{On the classification of constant mean 
curvature tori}, Ann. of Math \textbf{130} (1989), 
407--451.

\bibitem{RuhV}
E.~A.~Ruh and J.~Vilms,
\emph{The tension field of the Gauss map}, 
Trans. Amer. Math. Soc. \textbf{149} (1970), 
569--573.

\bibitem{Sch:cmclab}
N.~Schmitt, \emph{{CMCL}ab}, http://www.gang.umass.edu/software.

\bibitem{Sch:gallery}
\bysame, \emph{{CMC} {G}allery}, http://www.gang.umass.edu/cmcgallery.

\bibitem{Sch:tri}
\bysame, \emph{Constant mean curvature trinoids}, math.DG/0403036.

\bibitem{SteW}
I.~Sterling and H.~Wente, \emph{Existence and classification of constant mean
  curvature multibubbletons of finite and infinite type}, Indiana Univ. Math.
  J. \textbf{42} (1993), no.~4, 1239--1266.

\bibitem{TerU}
C.~Terng and K.~Uhlenbeck, \emph{B\"{a}cklund transformations and loop group
  actions}, Comm. Pure and Appl. Math \textbf{LIII} (2000), 1--75.

\bibitem{Uhl}
K.~Uhlenbeck, \emph{Harmonic maps into Lie groups 
  (Classical solutions of the chiral model)}, 
J. Diff. Geom. \textbf{30} (1989), 1--50.


\bibitem{Wu:dre}
H.~Wu, \emph{On the dressing action of loop groups on constant mean curvature
  surfaces}, T\^{o}hoku Math. J. \textbf{49} (1997), 599--621.

\bibitem{ZakS2} 
V.~E.~Zakharov and A.~S.~Shabat,
\emph{Integration of the nonlinear equations of 
mathematical physics by the method of the inverse 
scattering problem II}, 
Funktsional. Anal. i Prilozhen. \textbf{13} (1978),
13--22.

\end{thebibliography}
\end{document}